\documentclass[reqno,11pt]{amsart}
\usepackage{amsmath}
\usepackage{amsxtra}
\usepackage{amscd}
\usepackage{amsthm}
\usepackage{amsfonts}
\usepackage{amssymb}
\usepackage{eucal}
\usepackage[matrix,arrow,curve]{xy}
\theoremstyle{plain}
\newtheorem{Thm}[subsection]{Theorem}
\newtheorem{Cor}[subsection]{Corollary}
\newtheorem{Lem}[subsection]{Lemma}
\newtheorem{Prop}[subsection]{Proposition}
\newtheorem{Conj}[subsection]{Conjecture}

\theoremstyle{definition}
\newtheorem{Def}[subsection]{Definition}

\theoremstyle{remark}

\newtheorem{Rem}[subsection]{Remark}

\newtheorem{conj}{Conjecture}

\numberwithin{equation}{section}
\renewcommand{\rm}{\normalshape}

\newif\ifShowLabels
\ShowLabelstrue
\newdimen\theight
\def\TeXref#1{%
    \leavevmode\vadjust{\setbox0=\hbox{{\tt
        \quad\quad  {\small \rm #1}}}%
    \theight=\ht0
    \advance\theight by \lineskip
    \kern -\theight \vbox to
    \theight{\rightline{\rlap{\box0}}%
    \vss}%
    }}%

\ShowLabelsfalse% comment this out if labels should be printed

\renewcommand{\sec}[2]{\section{#2}\label{S:#1}%
    \ifShowLabels \TeXref{{S:#1}} \fi}
\newcommand{\ssec}[2]{\subsection{#2}\label{SS:#1}%
    \ifShowLabels \TeXref{{SS:#1}} \fi}

\newcommand{\refs}[1]{Section ~\ref{S:#1}}
\newcommand{\refss}[1]{Section ~\ref{SS:#1}}

\newcommand{\reft}[1]{Theorem ~\ref{T:#1}}
\newcommand{\refl}[1]{Lemma ~\ref{L:#1}}
\newcommand{\refp}[1]{Proposition ~\ref{P:#1}}

\newcommand{\refe}[1]{\eqref{E:#1}}
\newcommand{\refco}[1]{Conjecture ~\ref{Co:#1}}

\newenvironment{thm}[1]%
    { \begin{Thm} \label{T:#1}  \ifShowLabels \TeXref{T:#1} \fi }%
    { \end{Thm} }

\renewcommand{\th}[1]{\begin{thm}{#1} \sl }
\renewcommand{\eth}{\end{thm} }

\newenvironment{lemma}[1]%
    { \begin{Lem} \label{L:#1}  \ifShowLabels \TeXref{L:#1} \fi }%
    { \end{Lem} }
\newcommand{\lem}[1]{\begin{lemma}{#1} \sl}
\newcommand{\elem}{\end{lemma}}

\newenvironment{propos}[1]%
    { \begin{Prop} \label{P:#1}  \ifShowLabels \TeXref{P:#1} \fi }%
    { \end{Prop} }
\newcommand{\prop}[1]{\begin{propos}{#1}\sl }
\newcommand{\eprop}{\end{propos}}

\newenvironment{corol}[1]%
    { \begin{Cor} \label{C:#1}  \ifShowLabels \TeXref{C:#1} \fi }%
    { \end{Cor} }
\newcommand{\cor}[1]{\begin{corol}{#1} \sl }
\newcommand{\ecor}{\end{corol}}

\newenvironment{defeni}[1]%
    { \begin{Def} \label{D:#1}  \ifShowLabels \TeXref{D:#1} \fi }%
    { \end{Def} }
\newcommand{\defe}[1]{\begin{defeni}{#1} \sl }
\newcommand{\edefe}{\end{defeni}}

\newenvironment{remark}[1]%
    { \begin{Rem} \label{R:#1}  \ifShowLabels \TeXref{R:#1} \fi }%
    { \end{Rem} }
\newcommand{\rem}[1]{\begin{remark}{#1}}
\newcommand{\erem}{\end{remark}}

\newenvironment{conjec}[1]%
    { \begin{Conj} \label{Co:#1}  \ifShowLabels \TeXref{Co:#1} \fi }%
    { \end{Conj} }
\renewcommand{\conj}[1]{\begin{conjec}{#1} \sl }
\newcommand{\econj}{\end{conjec}}

\newcommand{\eq}[1]%
    { \ifShowLabels \TeXref{E:#1} \fi
       \begin{equation} \label{E:#1} }
\newcommand{\eeq}{ \end{equation} }

\newcommand{\prf}{ \begin{proof} }
\newcommand{\epr}{ \end{proof} }

\newcommand\alp{\alpha}     

\newcommand\gam{\gamma}

\newcommand\lam{\lambda}        \newcommand\Lam{\Lambda}

\newcommand\ome{\omega}     \newcommand\Ome{\Omega}

\newcommand\calM{{\mathcal{M}}}

\newcommand\calO{{\mathcal{O}}}

\newcommand\calW{{\mathcal{W}}}

\newcommand\PP{\mathbb{P}}

\newcommand\GG{\mathbb{G}}

\newcommand\CC{\mathbb{C}}

\newcommand\NN{\mathbb{N}}

 \newcommand\grg{{\mathfrak{g}}}

\newcommand\sdp{\times \hskip -0.3em {\raise 0.3ex
\hbox{$\scriptscriptstyle |$}}} % semidirect product

\newcommand\Gr{\operatorname{Gr}}

\newcommand\Hom{\operatorname {Hom}}

\newcommand\Int{\operatorname{Int}}

\newcommand\SL{{\rm SL}}

\newcommand\Sym{\operatorname{Sym}}

\newcommand\x{\times}
\newcommand\ten{\otimes}

\newcommand{\ra}{\rangle}
\newcommand{\la}{\langle}

\newcommand\nc{\newcommand}

\nc\aff{\operatorname{aff}}
\nc\oGr{\overline{\Gr}}
\nc\Bun{\operatorname{Bun}}
\nc\hgrg{\widehat{\grg}}
\renewcommand\Int{\operatorname{Int}}
\nc\bInt{\overline{\Int}}
\nc\hatLam{\widehat{\Lam}}
\nc\bmu{\overline{\mu}}
\nc\bnu{\overline{\nu}}
\nc\blambda{\overline{\lam}}
\renewcommand\SL{\operatorname{SL}}
\nc\ocalW{\overline{\calW}}
\nc\pos{\operatorname{pos}}
\nc\IH{\operatorname{IH}}
\nc\Rep{\operatorname{Rep}}
\nc\Gal{\operatorname{Gal}}
\nc{\tilGr}{\widetilde{\Gr}}

\nc\Pic{\operatorname{Pic}}
\nc{\HC}{{\mathcal{HC}}}
\nc{\on}{\operatorname}
\nc{\BA}{{\mathbb{A}}}
\nc{\BC}{{\mathbb{C}}}
\nc{\BF}{{\mathbb{F}}}
\nc{\BG}{{\mathbb{G}}}
\nc{\BM}{{\mathbb{M}}}
\nc{\BN}{{\mathbb{N}}}
\nc{\BQ}{{\mathbb{Q}}}
\nc{\BP}{{\mathbb{P}}}
\nc{\BR}{{\mathbb{R}}}
\nc{\BZ}{{\mathbb{Z}}}
\nc{\BS}{{\mathbb{S}}}

\nc{\CA}{{\mathcal{A}}}
\nc{\CB}{{\mathcal{B}}}
\nc{\CalC}{{\mathcal C}}
\nc{\CalD}{{\mathcal D}}
\nc{\CE}{{\mathcal{E}}}
\nc{\CF}{{\mathcal{F}}}
\nc{\CG}{{\mathcal{G}}}
\nc{\CH}{{\mathcal{H}}}
\nc{\CK}{{\mathcal{K}}}
\nc{\CL}{{\mathcal{L}}}
\nc{\CM}{{\mathcal{M}}}
\nc{\CMM}{{\mathcal{M}^{\operatorname{gen}}_\hbar(-\rho)}}
\nc{\CN}{{\mathcal{N}}}
\nc{\CO}{{\mathcal{O}}}
\nc{\CP}{{\mathcal{P}}}
\nc{\CQ}{{\mathcal{Q}}}
\nc{\CR}{{\mathcal{R}}}
\nc{\CS}{{\mathcal{S}}}
\nc{\CT}{{\mathcal{T}}}
\nc{\CU}{{\mathcal{U}}}
\nc{\CV}{{\mathcal{V}}}
\nc{\CW}{{\mathcal{W}}}
\nc{\CX}{{\mathcal{X}}}
\nc{\CY}{{\mathcal{Y}}}
\nc{\CZ}{{\mathcal{Z}}}

\nc{\gen}{{\operatorname{gen}}}
\nc{\cM}{{\check{\mathcal M}}{}}
\nc{\csM}{{\check{\mathcal A}}{}}
\nc{\obM}{{\overset{\circ}{\mathbf M}}{}}
\nc{\oCA}{{\overset{\circ}{\mathcal A}}{}}
\nc{\obA}{{\overset{\circ}{\mathbf A}}{}}
\nc{\ooM}{{\overset{\circ}{M}}{}}
\nc{\osM}{{\overset{\circ}{\mathsf M}}{}}
\nc{\vM}{{\overset{\bullet}{\mathcal M}}{}}
\nc{\nM}{{\underset{\bullet}{\mathcal M}}{}}
\nc{\obD}{{\overset{\circ}{\mathbf D}}{}}
\nc{\cp}{{\overset{\circ}{\mathbf p}}{}}
\nc{\ofZ}{{\overset{\circ}{\mathfrak Z}}{}}

\nc{\fa}{{\mathfrak{a}}}
\nc{\fb}{{\mathfrak{b}}}
\nc{\fg}{{\mathfrak{g}}}
\nc{\fgl}{{\mathfrak{gl}}}
\nc{\fh}{{\mathfrak{h}}}
\nc{\fj}{{\mathfrak{j}}}
\nc{\fm}{{\mathfrak{m}}}
\nc{\fn}{{\mathfrak{n}}}
\nc{\fu}{{\mathfrak{u}}}
\nc{\fp}{{\mathfrak{p}}}
\nc{\frr}{{\mathfrak{r}}}
\nc{\fs}{{\mathfrak{s}}}
\nc{\ft}{{\mathfrak{t}}}
\nc{\fT}{{\mathfrak{T}}}
\nc{\ofT}{{\overline{\mathfrak T}}}
\nc{\ofS}{{\overline{\mathfrak S}}}
\nc{\fsl}{{\mathfrak{sl}}}
\nc{\hsl}{{\widehat{\mathfrak{sl}}}}
\nc{\hgl}{{\widehat{\mathfrak{gl}}}}
\nc{\hg}{{\widehat{\mathfrak{g}}}}
\nc{\chg}{{\widehat{\mathfrak{g}}}{}^\vee}
\nc{\hn}{{\widehat{\mathfrak{n}}}}
\nc{\chn}{{\widehat{\mathfrak{n}}}{}^\vee}

\nc{\fA}{{\mathfrak{A}}}
\nc{\fB}{{\mathfrak{B}}}
\nc{\fD}{{\mathfrak{D}}}
\nc{\fE}{{\mathfrak{E}}}
\nc{\fF}{{\mathfrak{F}}}
\nc{\fG}{{\mathfrak{G}}}
\nc{\fI}{{\mathfrak{I}}}
\nc{\fJ}{{\mathfrak{J}}}
\nc{\fK}{{\mathfrak{K}}}
\nc{\fL}{{\mathfrak{L}}}
\nc{\fM}{{\mathfrak{M}}}
\nc{\fN}{{\mathfrak{N}}}
\nc{\frP}{{\mathfrak{P}}}
\nc{\fQ}{{\mathfrak Q}}
\nc{\fS}{{\mathfrak S}}
\nc{\fU}{{\mathfrak{U}}}
\nc{\fZ}{{\mathfrak{Z}}}

\nc{\bb}{{\mathbf{b}}}
\nc{\bc}{{\mathbf{c}}}
\nc{\be}{{\mathbf{e}}}
\nc{\bj}{{\mathbf{j}}}
\nc{\bn}{{\mathbf{n}}}
\nc{\bp}{{\mathbf{p}}}
\nc{\bq}{{\mathbf{q}}}
\nc{\bv}{{\mathbf{v}}}
\nc{\bx}{{\mathbf{x}}}
\nc{\by}{{\mathbf{y}}}
\nc{\bw}{{\mathbf{w}}}
\nc{\bA}{{\mathbf{A}}}
\nc{\bB}{{\mathbf{B}}}
\nc{\bC}{{\mathbf{C}}}
\nc{\bK}{{\mathbf{K}}}
\nc{\bD}{{\mathbf{D}}}
\nc{\bH}{{\mathbf{H}}}
\nc{\bM}{{\mathbf{M}}}
\nc{\bN}{{\mathbf{N}}}
\nc{\bO}{{\mathbf{O}}}
\nc{\bQ}{{\mathbf Q}}
\nc{\bS}{{\mathbf{S}}}
\nc{\bT}{{\mathbf{T}}}
\nc{\bV}{{\mathbf{V}}}
\nc{\bW}{{\mathbf{W}}}
\nc{\bX}{{\mathbf{X}}}
\nc{\bP}{{\mathbf{P}}}
\nc{\bZ}{{\mathbf{Z}}}

\nc{\sA}{{\mathsf{A}}}
\nc{\sB}{{\mathsf{B}}}
\nc{\sC}{{\mathsf{C}}}
\nc{\sD}{{\mathsf{D}}}
\nc{\sF}{{\mathsf{F}}}
\nc{\sK}{{\mathsf{K}}}
\nc{\sM}{{\mathsf{M}}}
\nc{\sO}{{\mathsf{O}}}
\nc{\sQ}{{\mathsf{Q}}}
\nc{\sP}{{\mathsf{P}}}
\nc{\sV}{{\mathsf{V}}}
\nc{\sW}{{\mathsf{W}}}
\nc{\sZ}{{\mathsf{Z}}}
\nc{\sfp}{{\mathsf{p}}}
\nc{\sr}{{\mathsf{r}}}
\nc{\st}{{\mathsf{t}}}
\nc{\sx}{{\mathsf{x}}}
\nc{\sfb}{{\mathsf{b}}}
\nc{\sfc}{{\mathsf{c}}}
\nc{\sd}{{\mathsf{d}}}
\nc{\sg}{{\mathsf{g}}}
\nc{\sfl}{{\mathsf{l}}}

\nc{\BK}{{\bar{K}}}
\nc{\balpha}{{\bar{\alpha}}}

\nc{\tA}{{\widetilde{\mathbf{A}}}}
\nc{\tB}{{\widetilde{\mathcal{B}}}}
\nc{\tg}{{\widetilde{\mathfrak{g}}}}
\nc{\tG}{{\widetilde{G}}}
\nc{\TM}{{\widetilde{\mathbb{M}}}{}}
\nc{\tO}{{\widetilde{\mathsf{O}}}{}}
\nc{\tU}{{\widetilde{\mathfrak{U}}}{}}
\nc{\TZ}{{\tilde{Z}}}
\nc{\tZ}{\widetilde{Z}{}}
\nc{\tx}{{\tilde{x}}}
\nc{\tbv}{{\tilde{\bv}}}
\nc{\tfP}{{\widetilde{\mathfrak{P}}}{}}
\nc{\tz}{{\tilde{\zeta}}}
\nc{\tmu}{{\tilde{\mu}}}

\nc{\td}{\ddot{\underline{d}}{}}
\nc{\tzeta}{\widetilde{\zeta}{}}
\nc{\hd}{{\widehat{\underline{d}}}}
\nc{\hG}{{\widehat{G}}}
\nc{\hBP}{\widehat{\mathbb P}{}}
\nc{\hQ}{{\widehat{Q}}}
\nc{\hsM}{\widehat{\mathsf M}{}}
\nc{\hfM}{\widehat{\mathfrak M}{}}
\nc{\hCP}{\widehat{\mathcal P}{}}
\nc{\hCR}{\widehat{\mathcal R}{}}
\nc{\hCS}{{\widehat{\mathcal S}}}
\nc{\hfZ}{\widehat{\mathfrak Z}{}}

\nc{\urho}{\underline{\rho}}
\nc{\uB}{\underline{B}}
\nc{\uC}{{\underline{\mathbb{C}}}}
\nc{\ui}{\underline{i}}
\nc{\ofP}{{\overline{\mathfrak{P}}}}
\nc{\hrho}{{\hat{\rho}}}

\nc{\unl}{\underline}
\nc{\ol}{\overline}
\nc{\one}{{\mathbf{1}}}
\nc{\two}{{\mathbf{t}}}

\nc{\Tot}{{\mathop{\operatorname{\rm Tot}}}}
\nc{\Hilb}{{\mathop{\operatorname{\rm Hilb}}}}
\nc{\CHom}{{\mathop{\operatorname{{\mathcal{H}}\it om}}}}
\nc{\defi}{{\mathop{\operatorname{\rm def}}}}
\nc{\length}{{\mathop{\operatorname{\rm length}}}}
\nc{\Cliff}{{\mathsf{Cliff}}}
\nc{\Fl}{{\mathsf{Fl}}}
\nc{\Fib}{{\mathsf{Fib}}}
\nc{\Coh}{{\mathsf{Coh}}}
\nc{\FCoh}{{\mathsf{FCoh}}}

\nc{\reg}{{\text{\rm reg}}}

\nc{\cplus}{{\mathbf{C}_+}}
\nc{\cminus}{{\mathbf{C}_-}}
\nc{\cthree}{{\mathbf{C}_*}}
\nc{\Qbar}{{\bar{Q}}}

\nc{\bh}{{\bar{h}}}
\nc{\bOmega}{{\overline{\Omega}}}
\nc\tGr{\widetilde{\Gr}}

\nc{\seq}[1]{\stackrel{#1}{\sim}}
\nc\ogu{\overline{G/U}}
\nc\chlam{{\check{\lam}}}
\nc\St{\operatorname{St}}

\nc\uS{\underline{S}}
\nc\QM{\mathcal{QM}}
\nc{\chmu}{\check{\mu}}
\begin{document}
\title{Macdonald polynomials, Laumon spaces and perverse coherent sheaves}
\author{Alexander Braverman, Michael Finkelberg and Jun'ichi Shiraishi}

\begin{abstract}
Let $G$ be an almost simple simply connected complex Lie group,
and let $G/U_-$ be its base affine space.
In this paper we formulate a conjecture, which provides a new geometric interpretation
of the Macdonald polynomials associated to $G$ via perverse coherent sheaves on the scheme of formal arcs in the affinization of $G/U_-$.
%so called {\it spaces of quasi-maps} from $\PP^1$ to $\CB$.
We prove our conjecture for $G=SL(N)$
using the so called Laumon resolution of the space of quasi-maps (using this resolution one can reformulate the statement so that only ``usual" (not perverse) coherent sheaves are used).
In the course of the proof we also give a $K$-theoretic version of the main
result of \cite{Negut}.
\end{abstract}
\maketitle

\sec{int}{Introduction}

\ssec{}{Notations}Let $\grg$ be a semi-simple Lie algebra over $\CC$ and let $G$ be the corresponding simply connected
group. Let $B,B_-\subset G$ be a pair of opposite Borel subgroups with unipotent radicals $U,U_-$
and let $T= B\cap B^-$ be the corresponding maximal torus.
We denote by $\Lam$ the lattice of cocharacters of $T$ (this is also the coroot lattice of $G$, since $G$ is simply
connected) and by $\check\Lam$ the lattice of characters of $T$. We denote by $\Lam_+$ the cone consisting of sums
of positive coroots of $G$ with non-negative coefficients. Similarly, we denote by $\check\Lam^+$ the cone of dominant
weights.

We denote by $\CB$ the flag variety of $G$. It can be identified with the quotient $G/B$. The choice of $B_-$ gives a point
in the  open $B$-orbit in $\CB$.

For a pair of variables $p,q$ and for any $n\in \NN\cup \infty$ we set
$$
(p;q)_n:=(1-p)(1-qp)\ldots(1-q^{n-1}p).
$$

%---------------------------------------------------------------
\ssec{}{Quasi-maps and Laumon spaces}For $\alp\in\Lam_+$ we denote by $_{\fg}\calM^{\alp}$ the moduli space
of maps $\PP^1\to \CB$ of degree $\alp$ and by $_{\fg}\QM^{\alp}$ its quasi-maps compactification (cf. \cite{Bicm} for
a survery on quasi-maps); we shall sometimes omit the subscript $\fg$ when it does not lead to a confusion.
The scheme $\QM^{\alp}$ possesses a natural stratification
$$
\QM^{\alp}=\bigsqcup\limits_{0\leq\beta\leq \alp}\calM^{\beta}\x\Sym^{\alp-\beta}(\PP^1),
$$
where $\Sym^{\alp-\beta}(\PP^1)$ stands for the space of all formal linear combinations $\sum\gam_i x_i$
where $\gam_i\in \Lam_+$, $x_i\in \PP^1$ and $\sum \gam_i=\alp$. The points $\{ x_i\}$ are called the
points of defect of the corresponding quasi-map.

Similarly, we denote by $Z^{\alp}$ the space of based quasi-maps of degree $\alp$ (i.e. those quasi-maps, which have no defect
at $\infty\in \PP^1$ and which send $\infty$ to $B_-$ regarded as a point in $\CB$).
The space $\QM^{\alp}$ has a natural action of $PGL(2)\x G$; here the first factor acts on $\PP^1$ and the second
on $\CB$. This action does not preserve $Z^{\alp}$; however, $\GG_m\x T$ still acts on $Z^{\alp}$.

It is well-known that the space $\QM^{\alp}$ is usually singular, but when $G=\SL(N)$ it has a natural small
resolution of singularities by means of Laumon's quasiflags' space $\CQ^{\alp}$. By the definition, it consists of flags
$$
0\subset \calW_1\subset\calW_2\subset \cdots\subset \calW_N=\calO_{\PP^1}^N,
$$
where $\calW_i$ is a locally free sheaf on $\PP^1$ of rank $i$ and such that
$$
\deg\calW_i=-\la \alp,\check\omega_i\ra.
$$
We shall denote by $\fQ^{\alp}$ the corresponding ``based" version of $\CQ^{\alp}$.

As before, $\CQ^{\alp}$ has a natural action of $PGL(2)\x G$ and $\fQ^{\alp}$ has a natural action of $\GG_m\x T$.
%------------------------------------------------------------------------------------------------
\ssec{}{Geometric interpretation of the ``Macdonald function" for $G=\SL(N)$}
In the case $G=\SL(N)$ we identify $\Lam_+$ with $\NN^{N-1}$ by using the
simple coroots $\alp_i$ as a basis of $\Lam$.
Similarly, we identify $\check\Lam^+$ with $\NN^{N-1}$ by using the fundamental
weights $\check\omega_i$ as a basis. Also we have the natural isomorphism
$T\simeq \GG_m^{N-1}$.

For any $\alp\in\Lam_+$ let us set
\eq{jalpha}\fJ_\alpha(q,t,z)=[H^\bullet(\fQ^\alpha,\Omega^\bullet_{\fQ^\alpha})]:=
\sum_{i,j}(-1)^{i+j}t^j[H^i(\fQ^\alpha,\Omega^j_{\fQ^\alpha})].
\end{equation}
Here $[H^i(\fQ^\alpha,\Omega^j_{\fQ^\alpha})]$ means the character of $H^i(\fQ^\alpha,\Omega^j_{\fQ^\alpha})$
as a representation of $\GG_m\x T$; in other words, it is a function of $q\in \GG_m$ and $z\in T$. More precisely, the coordinate functions $z_i,\ i=1,\ldots,
N-1$, satisfy $\check\omega_i=z_1\cdots z_i$.

We would like to organize all the $\fJ_{\alp}$ into a generating function. Namely, let us set:
$$
J(q,t,z,x)=\sum_{\alpha\in\BN^{N-1}}x^\alpha\fJ_\alpha(q,t,z);\qquad
\fJ(q,t,z,x)=\prod_{i=1}^{N-1}x_i^{\log(\check\omega_i)/\log q}J(q,t,z,x).
$$
Also, for $1\leq i\leq N$, we consider the difference operator $T_{i,q^{\pm1}}$
%on functions of $q,t,z,x$
defined as follows: $T_{i,q^{\pm1}}F(q,t,z,x_1,\ldots,x_{N-1}):=
F(q,t,z,x_1,\ldots,x_{i-2},q^{\mp1}x_{i-1},q^{\pm1}x_i,x_{i+1},\ldots,x_{N-1})$.
Our first main result is the following
\th{shir-int}
\begin{enumerate}
\item
Define the function $z_N$ on the Cartan torus $T$ of $\SL(N)$ by
$z_N:=z_1^{-1}\cdots z_{N-1}^{-1}$. Then we have
$$
\sD\fJ(q,t,z,x)=(z_1+\ldots+z_N)\fJ(q,t,z,x),
$$
where
$$
\sD:=\sum_{i=1}^N\prod_{j<i}
\frac{1-q^{-1}t^{i-j-1}x_j\cdots x_{i-1}}{1-t^{i-j}x_j\cdots x_{i-1}}
\prod_{k>i}\frac{1-qt^{k-i+1}x_i\cdots x_{k-1}}{1-t^{k-i}x_i\cdots x_{k-1}}
T_{i,q^{-1}}
$$
\item
$$
\lim_{\alpha\to\infty}\fJ_\alpha(q,t,z)=\prod_{1\leq i<j\leq N}
{(qt z_j/z_i;q)_{\infty } \over (q z_j/z_i;q)_{\infty}}\times
\left({(qt;q)_{\infty } \over (q;q)_{\infty}}\right)^{N-1}\times
\prod_{i=1}^{N-2}
\left({(qt^{i+1};q)_\infty \over (t^i;q)_\infty}\right)^{N-i-1}.
$$
\end{enumerate}
\eth
Some remarks about \reft{shir-int} are in order. First, the operator $\sD$ is a version of one of the
Macdonald difference operators; it is easy to see that the first assertion of \reft{shir-int} implies that
$\fJ$ is an eigen-function of all the (suitably normalized) Macdonald
operators and thus (up to some normalization factor) it is
equal to the {\em Baker-Akhiezer} function for the Macdonald operators in the terminology of \cite{ES} or \cite{CE}; it is
also often called {\em the Macdonald function}.
Moreover, the second assertion can be deduced from the first one and the
results of \cite{ES}, \cite{CE}, but we are going to
give an independent proof of this result.

It should also be noted that some limiting cases of \reft{shir-int} have been known before. In particular,
the case $t=0$ is treated in \cite{BF0} (cf. also
\cite{BF11} for a generalization to arbitrary $G$).
Also, in \cite{Negut} the $q\to 1$ version of \reft{shir-int}
is proved. It should be noted that the proofs in {\em loc. cit.} are representation-theoretic: they are based on an interpretation of the (localized) equivariant $K$-theory (resp.
localized equivariant cohomology) of all the $\fQ^{\alp}$
as the universal Verma module for the quantum group $U_q({\mathfrak{sl}}(N))$
(resp. of the lie algebra ${\mathfrak{sl}}(N)$). On the other hand, the proof
of \reft{shir-int} given in this paper
is purely computational: using Atiyah-Bott-Lefschetz localization
formula one can produce a combinatorial expression for the
function $\fJ_{\alp}$ and thus reduce \reft{shir-int}(1) to
a combinatorial identity, which can be proven by an explicit (but fairly long) computation. It would be very interesting to extend the methods of {\em loc. cit.} to the present situation.
%----------------------------------------------------------------------------------------------------------------
\ssec{}{Geometric interpretation of Macdonald polynomials for $G=\SL(N)$}
The Macdonald operators are usually used in order to define the so called Macdonald polynomials. This is a series
of $W$-invariant polynomials $P_{\check \lam}(q,t,z)$ on the torus $T$ (recall that $z\in T$)
depending on a dominant weight $\check\lam\in\check\Lam^+$ and on the variables $q,t\in \GG_m$.
We would like to present a geometric construction of these polynomials. Let us explain how to
do it in the $\SL(N)$-case. The conjectural generalization to arbitrary $G$ is discussed in the next Subsection.

First, for any $\check\lam\in\check\Lam$ one can construct a line bundle $\calO(\check\lam)$ on $\QM^{\alp}$;
abusing the notation we are going to denote its pull-back to $\CQ^{\alp}$ also by $\calO(\check \lam)$.
The construction is discussed in~\cite{BF11}.
We are not going to recall the construction in the Introduction, but let us just note that it requires a choice of
a point $\infty\in\PP^1$. Hence, the bundle $\calO(\check\lam)$ is not $PGL(2)$-equivariant. However, it is still
equivariant with respect to the diagonal torus $\GG_m\subset PGL(2)$. In particular, it makes sense to consider the character of
$H^\bullet(\CQ^{\alp}, \Ome^\bullet_{\CQ^{\alp}}\ten\calO(\check\lam))$ with respect to the action of $\GG_m\x G$, which we shall denote by $[H^\bullet(\CQ^{\alp}, \Ome^\bullet_{\CQ^{\alp}}\ten\calO(\check\lam))]$. By definition this character is
$W$-invariant function on $\GG_m\x T$.
%------------------------------------------------------------
\th{mac-int}
\begin{enumerate}
\item
Assume that $\check\lam\in\check\Lam$ is not dominant. Fix
$j,k\in \NN$.
Then for $\alp$ sufficiently large we have
$$
H^k(\CQ^{\alp},\Ome^j_{\CQ^{\alp}}\ten \calO(\check\lam))=0.
$$
\item
For any $\check\lam\in\check\Lam$ there exists the limit
$\lim\limits_{\alp\to \infty}[H^\bullet(\CQ^{\alp}, \Ome^\bullet_{\CQ^{\alp}}\ten\calO(\check\lam))]$.
We shall denote the above limit by $H_{\check\lam}(q,t,z)$.
Note that it follows from the first assertion that $H_{\check\lam}=0$ when $\check\lam$ is not dominant.
\item
$$
H_0(q,t,z)=\frac{(1+t)(1+t+t^2)\ldots(1+t+\ldots+t^{N-1})}
{(1-t^{N-1})^2(1-t^{N-2})^4\ldots(1-t^3)^{2N-6}(1-t^2)^{2N-4}}\cdot
\frac{1}{(1-t^N)(1-t)^{N-2}}.
$$
\item
For any $\chlam=\sum l_i\check\ome_i\in\check\Lam^+$
(here $\check\ome_i$ denotes the $i$-th fundamental weight of
$\SL(N)$) we have
$$
H_\chlam=H_0\prod_{1\leq i\leq j\leq N-1}\frac{(t^{j-i+1};q)_{l_i+\ldots+l_j}}
{(t^{j-i}q;q)_{l_i+\ldots+l_j}}~P_\chlam.
$$
In other words, $H_\chlam$ is equal to $P_\chlam$ up to an explicit factor.
\end{enumerate}
\eth

\ssec{}{The case of arbitrary $G$}
In this subsection we are going to give a conjectural
formulation\footnote{The reader should be warned that we do not know how
to formulate a version of \reft{shir-int} for arbitrary $G$.}
of \reft{mac-int} for arbitrary $G$.
The formulation is based on the theory of perverse coherent sheaves developed by D.~Arinkin and R.~Bezrukavnikov
(cf.~\cite{AB}). For simplicity, in this Introduction we shall assume that $G$ is simply laced (in the general case certain modification of the construction given below is needed; the details are explained in \refs{spec}).

First let us introduce the infinite type scheme $_\fg\bQ$ (discussed also
in~\cite[Section~2.2]{BF12}): it is
the quotient by the action of the Cartan torus $T\subset G$ of the space of
maps from $\on{Spec}R=\on{Spec}\BC[[{\mathbf t}^{-1}]]$ to the affinization
of the base affine space $\overline{G/U_-}$ taking value in $G/U_-$ at the
generic point. This scheme is equipped with the action of the proalgebraic group
$G(R)$; the open orbit $_\fg\bQ_\infty=\ {}_\fg\bQ^0$ is nothing but
$G(R)/T\cdot U_-(R)$: the maps taking value in $G/U_-$ at the closed point
$r\in\on{Spec}R$. We denote by $\fj$ the open embedding of $_\fg\bQ^0$
into $_\fg\bQ$.
All the $G(R)$-orbits in $_\fg\bQ$ are numbered by the
defects at $r$ taking value in the cone of positive coroots
$\Lambda_+$ of $G:\ {}_\fg\bQ=\bigsqcup_{\alpha\in\Lambda_+}\ {}_\fg\bQ^\alpha$.
The codimension of $_\fg\bQ^\alpha$ in $_\fg\bQ$ equals $2|\alpha|$.

We introduce the perversity $p(\ {}_\fg\bQ^\alpha)=|\alpha|$; it is immediate
that the function $p$ is strictly monotone and comonotone in the sense
of~\cite{AB}. For a locally free $G(R)\rtimes\BG_m$-equivariant sheaf
$\CF$ on $_\fg\bQ^0$ the construction of~\cite[Section~4]{AB} produces an
object $\fj_{!*}\CF$ of $G(R)\rtimes\BG_m$-equivariant quasicoherent derived
category on $_\fg\bQ$.

\conj{76-int}
(a) For a nondominant $G$-weight $\chlam$ we have
$[H^\bullet(_\fg\bQ,\fj_{!*}
(\Omega^\bullet_{_\fg\bQ{}^0})
\otimes\CO(\chlam))]=0$.

(b) For a dominant $G$-weight $\chlam$ we have
$$
[H^\bullet(_\fg\bQ,\fj_{!*}
(\Omega^\bullet_{_\fg\bQ{}^0})
\otimes\CO(\chlam))]=H_0\prod_{\alpha\in R^+(\check\fg)}
\frac{(t^{|\alpha|};q)_{\langle\alpha,\chlam\rangle}}
{(t^{|\alpha|-1}q;q)_{\langle\alpha,\chlam\rangle}}
\prod\frac{(t^{|\alpha|-1};q)_\infty}{(qt^{|\alpha|};q)_\infty}
P_\chlam
$$
where $P_\chlam(q,t,z)$ is the Macdonald polynomial for $G$,
and the second product is taken over all {\em nonsimple} positive roots
of $R^+(\check\fg)$.
\econj
We explain in \refs{spec} why \refco{76-int} is equivalent
to \reft{mac-int} for $G=\SL(N)$.

%-------------------------------------------------------------

\ssec{}{Organization of the paper}In \refs{comnot} and \refs{macf} we gather some combinatorial information about
Macdonald polynomials and the ``Macdonald function" for root systems of type A. In \refs{van} we prove a generalization
of the Sommese vanishing theorem, which in particular implies \reft{mac-int}(1). In \refs{DR} and \refs{DE} we prove~\reft{shir-int} and~\reft{mac-int}.
Finally, in \refs{spec} we give a careful formulation of \refco{76-int} for arbitrary $G$ and show that for $G=\SL(N)$ it is equivalent to \reft{mac-int}.
%--------------------------------------------------------------------------------------------------------------
\ssec{ack}{Acknowledgments}
We are grateful to P.~Etingof and B.~Feigin for very useful discussions and
introducing the second and the third authors to each other. The vanishing
theorems of~\refs{van} are the results of generous explanations by
E.~Amerik, M.~Brion, D.~Kaledin, and especially S.~Kov\'acs.
The speculations of~\refs{spec} are due to patient explanations by
D.~Arinkin and R.~Bezrukavnikov. J.S. is very grateful to M.~Noumi for
stimulating discussion and collaboration.

M.~F. was partially supported by the RFBR grants 12-01-33101, 12-01-00944, the National Research University Higher School of Economics' Academic Fund award No.12-09-0062 and
the AG Laboratory HSE, RF government grant, ag. 11.G34.31.0023.
This study was carried out within the National Research University Higher School of Economics
Academic Fund Program in 2012-2013, research grant No. 11-01-0017.
This study comprises research findings from the ``Representation Theory
in Geometry and in Mathematical Physics" carried out within The
National Research University Higher School of Economics' Academic Fund Program
in 2012, grant No 12-05-0014.
Research of J.S. is supported by the Grant-in-Aid for Scientific
Research C-24540206.

%----------------------------------------------------------------------------------------------------------

\sec{comnot}{Combinatorial notations}

\ssec{mapo}{Macdonald polynomials}
We follow the notations in~\cite{Mac} (especially, part VI), cf. also~\cite{NS}.
Let $N$ be a positive integer and $q,\st$ be independent indeterminates.
Let $\Lambda_{N,\BF}$ be the ring of symmetric polynomials in $N$ variables
with coefficients in $\BF={\BQ}(q,\st)$. Set
\begin{align}
T_{q,y_i}f(y_1,\ldots,y_N)=f(y_1,\ldots,q y_i,\ldots,y_N).
\end{align}

For a partition $\lambda$, the Macdonald polynomial $P_\lambda(y;q,\st)\in \Lambda_{N,F}$
is uniquely characterized by the conditions:
\begin{align}
&P_\lambda=m_\lambda+\sum_{\mu<\lambda}u_{\lambda\mu} m_\mu,\\
&\sD^1_N P_\lambda=\sum_{i=1}^N q^{\lambda_i}\st^{N-i} \cdot  P_\lambda,
\end{align}
where $m_\lambda$ is the monomial symmetric function, and
$\sD^1_N=\sD^1_N(q,\st) $ is the Macdonald difference operator
\begin{align}
\sD^1_N=\sum_{i=1}^N \prod_{j\neq i}{\st y_i-y_j\over y_i-y_j} T_{q,y_i}.
\end{align}

%Let $m,n$ be positive integers, and
%$y=(y_1,\ldots,y_m),y=(y_1,\ldots,y_n)$ be sets of variables.
%We have the kernel identity of Cauchy type
%\begin{align}
%&\Pi(y,y;q,\st)=\prod_{i,j}{(\st y_iy_j;q)_\infty \over (y_iy_j;q)_\infty},\\
%&\left(\st^{-m}\sD^1_{m,y}-\sum_{i=1}^m \st^{-i}\right)\Pi(y,y)
%=
%\left(\st^{-n}\sD^1_{n,y}-\sum_{i=1}^n \st^{-i}\right)\Pi(y,y).
%\end{align}

%%%%%%%%%%%%%%%%%%%%%%%%%%%%%%%%%%%%%%%%%%%%%%%%%
\ssec{tab}{Tableau}
Let
$\lambda=(\lambda_1,\lambda_2,\cdots),\mu=(\mu_1,\mu_2,\cdots)$
be partitions satisfying $\mu\subset \lambda$.
The necessary and sufficient condition for the
skew diagram
$\theta=\lambda-\mu$
to be a horizontal strip is
\begin{align}
\lambda_1\geq \mu_1\geq\lambda_2\geq \mu_2\geq\cdots.
\end{align}
This can be written as
\begin{align}
0\leq \lambda_i-\mu_i\leq \lambda_i-\lambda_{i+1}\qquad (i\geq 1).
\end{align}

A (column-strict) tableau $T$ of shape $\lambda$ is defined to be a sequence of partitions
\begin{align}
 \phi=\lambda^{(0)}\subset \lambda^{(1)}\subset
\cdots \lambda^{(N)}=\lambda\label{tableau}
\end{align}
such that every skew diagram $\theta^{(i)}=\lambda^{(i)}- \lambda^{(i-1)}$
is a horizontal strip.
Writing $\lambda^{(i)}=(\lambda^{(i)}_1,\lambda^{(i)}_2,\cdots)$,
the condition for the $T$ being a tableau reads
\begin{align}
0\leq \lambda_i^{(j)}-\lambda_i^{(j-1)}\leq \lambda_i^{(j)}-\lambda_{i+1}^{(j)}\qquad
 (1\leq i,1\leq j\leq N)\label{futou-1}.
\end{align}
Note that from $ \lambda^{(0)}=\phi$ and the inequality (\ref{futou-1}) we have
\begin{align}
\lambda_i^{(j)}=0
\qquad  (i>j)\label{kousoku-1}.
\end{align}

For each skew diagram $\theta^{(i)}=\lambda^{(i)}- \lambda^{(i-1)}$, set
\begin{align}
&
\theta_{i,j}=\lambda^{(j)}_i- \lambda^{(j-1)}_i\qquad (1\leq i\leq N,1\le j\leq N),
\end{align}
for simplicity of display. Then the constraint (\ref{kousoku-1}) means
\begin{align}
&\theta_{i,j}=0 \qquad (i>j),\label{kousoku-2}\\
&\lambda_i=\sum_{k=i}^N \theta_{i,k}\qquad(1\leq i\leq N).\label{kousoku-3}
\end{align}

Hence the tableau $T$ uniquely gives us a set of
$N(N-1)/2$
nonnegative integers $\{\theta_{i,j}|1\leq i<j\leq N\}$ satisfying (\ref{futou-1}), namely
\begin{align}
0\leq \theta_{i,j}\leq
\lambda_{i}-\lambda_{i+1}-
\sum_{k=j+1}^{N} \left(\theta_{i,k}-\theta_{i+1,k}\right)
\qquad(1\leq i<j\leq N)\label{futou-2}.
\end{align}
Conversely, a set of nonnegative integers $\{\theta_{i,j}\}$
satisfying (\ref{futou-2}) uniquely
gives us a sequence of partitions $\lambda^{(j)}=(\lambda^{(j)}_1,\lambda^{(j)}_2,\ldots)$
\begin{align}
\lambda_i^{(j)}=\sum_{k=1}^j \theta_{i,k},
\end{align}
which is a tableau.

It is convenient to consider a set of $N\times N$ upper triangular matrices $\mathsf{M}^{(N)}$
having $\{\theta_{i,j}\}$'s as nonzero entries, and zeros on the diagonal:
\begin{align}
\label{215}
\mathsf{M}^{(N)}=\{\theta=(\theta_{i,j})_{1\leq i,j\leq N}|
\theta_{i,j}\in {{\mathbb Z}_{\geq 0}},\theta_{i,j}=0 \mbox{ if }i\geq j\}.
\end{align}
We have a natural projection $\sM^{(N)}\to\sM^{(N-1)}$ forgetting the last column.

\lem{sh1}
Let
$\lambda=(\lambda_1,\ldots,\lambda_N)$ be a partition.
We have a one to one mapping from
the set of (column-strict)  tableaux of
shape $\lambda$ to the elements in the
polyhedral region ${\mathsf{Pol}}_\lambda\in\mathsf{M}^{(N)}$
defined by
\begin{align}
{\mathsf{Pol}}_\lambda=\{\theta\in\mathsf{M}^{(N)}|
0\leq \theta_{i,j}\leq
\lambda_{i}-\lambda_{i+1}-
\sum_{k=j+1}^{N} \left(\theta_{i,k}-\theta_{i+1,k}\right)\}.
\end{align}
\elem

\lem{sh2}
The size of the
skew diagram $\theta^{(i)}=\lambda^{(i)}- \lambda^{(i-1)}$ is written as
\begin{align}
|\theta^{(i)}|=\lambda_i+\sum_{a=1}^{i-1}\theta_{a,i}-
\sum_{b=i+1}^{N}\theta_{i,b}.
\end{align}
\elem

%%%%%%%%%%%%%%%%%%%%%%%%%%%%%%%%%%%%%%%%%%%%%%
\ssec{tsf}{Tableaux sum formula}
We recall the tableaux sum formula for the Macdonald polynomials.

The Macdonald polynomial $P_\lambda$ is written as
\begin{align}
P_\lambda=\sum_T \psi_T(q,\st) y^T.
\end{align}
where $T$ runs over the set of tableaux of shape $\lambda$,
$y^T$ denotes the monomial defined in terms of the
weights $\alpha=(|\theta^{(1)}|,|\theta^{(2)}|,\ldots,|\theta^{(N)}|)$ of $T$ as
\begin{align}
y^T=y^\alpha=y^\lambda \prod_{1\leq i<j\leq N}(y_j/y_i)^{\theta_{i,j}},
\end{align}
and the coefficient $\psi_T(q,\st)$ is given by
\begin{align}
&\psi_T(q,\st)=\prod_{i=1}^N \psi_{\lambda^{(i)}/ \lambda^{(i-1)}}(q,\st),\\
&\psi_{\lambda/\mu}=\prod_{1\leq i\leq j\leq \ell(\mu)}
{f(q^{\mu_i-\mu_j}\st^{j-i})f(q^{\lambda_i-\lambda_{j+1}}\st^{j-i}) \over
f(q^{\lambda_i-\mu_j}\st^{j-i})f(q^{\mu_i-\lambda_{j+1}}\st^{j-i}) },\\
&f(u)={(\st u;q)_\infty\over (qu;q)_\infty}.
\end{align}

For a nonnegative integer $\theta\in {{\mathbb Z}_{\geq 0}}$, we have
\begin{align*}
%&&{f(u)\over f(q^\theta u)}={(\st u;q)_\theta\over (qu;q)_\theta},\\
&&{f(u)\over f(q^{-\theta} u)}={(q^{-\theta+1}u;q)_\theta\over (q^{-\theta}\st u;q)_\theta}=
(q/\st)^\theta{(1/u;q)_\theta\over (q/\st u;q)_\theta}.
\end{align*}
where $(p;q)_n:=(1-p)(1-qp)\ldots(1-q^{n-1}p)$.
Hence we have
\begin{align}
&\psi_T(q,\st)\nonumber \\
&=\prod_{k=1}^{N}
\prod_{1\leq i\leq j\leq k-1}
{f(q^{\lambda^{(k-1)}_i-\lambda^{(k-1)}_j}\st^{j-i})
f(q^{\lambda^{(k)}_i-\lambda^{(k)}_{j+1}}\st^{j-i}) \over
f(q^{\lambda^{(k)}_i-\lambda^{(k-1)}_j}\st^{j-i})
f(q^{\lambda^{(k-1)}_i-\lambda^{(k)}_{j+1}}\st^{j-i}) }\\
&=\prod_{k=1}^{N}
\prod_{1\leq i\leq j\leq k-1}
{f(q^{-\theta_{i,k}+\lambda^{(k)}_i-\lambda^{(k-1)}_j}\st^{j-i})
f(q^{\lambda^{(k)}_i-\lambda^{(k)}_{j+1}}\st^{j-i}) \over
f(q^{\lambda^{(k)}_i-\lambda^{(k-1)}_j}\st^{j-i})
f(q^{-\theta_{i,k}+\lambda^{(k)}_i-\lambda^{(k)}_{j+1}}\st^{j-i}) }\nonumber\\
&=
\prod_{k=1}^{N}
\prod_{1\leq i\leq j\leq k-1}
{(q^{-\lambda^{(k)}_i+\lambda^{(k-1)}_j+1}\st^{-j+i-1};q)_{\theta_{i,k}}\over
(q^{-\lambda^{(k)}_i+\lambda^{(k-1)}_j}\st^{-j+i};q)_{\theta_{i,k}}}
{(q^{-\lambda^{(k)}_i+\lambda^{(k)}_{j+1}}\st^{-j+i};q)_{\theta_{i,k}}\over
(q^{-\lambda^{(k)}_i+\lambda^{(k)}_{j+1}+1}\st^{-j+i-1};q)_{\theta_{i,k}}}.\nonumber
\end{align}

%%%%%%%%%%%%%%%%%%%%%%%%%%%%%%%%%%%%%%%%%%%%%%%%%
%%%%%%%%%%%%%%%%%%%%%%%%%%%%%%%%%%%%%%%%%%%%%%%%%
\sec{macf}{Macdonald function}
%%%%%%%%%%%%%%%%%%%%%%%%%%%%%%%%%%%%%%%%%%%%%%%%%
\ssec{mhts}{Multiple hypergeometric-type series}
Let
$q,\st,z_1,z_2,\ldots,z_N$ be independent indeterminates.
Recall the projection $\sM^{(N)}\to\sM^{(N-1)}$, see the line after~(\ref{215}).
Define a sequence of rational functions
$c_N(\theta;z_1,\ldots,z_N;q,\st)\in {\mathbb Q}(q,\st,z_1,\ldots,z_N)$
inductively as follows:
\begin{align}
&c_1(-;z_1;q,\st)=1,\\
&c_N(\theta\in \mathsf{M}^{(N)};z_1,\ldots,z_N;q,\st)\nonumber\\
&=c_{N-1}(\theta\in \mathsf{M}^{(N-1)};
q^{-\theta_{1,N}}z_1,\ldots,q^{-\theta_{N-1,N}}z_{N-1};q,\st)\\
&\times
\prod_{1\leq i\leq j\leq N-1}
{(\st z_{j+1}/z_i;q)_{\theta_{i,N}}\over (q z_{j+1}/z_i;q)_{\theta_{i,N}}}
{(q^{-\theta_{j,N}}q z_j/\st z_i;q)_{\theta_{i,N}}\over (q^{-\theta_{j,N}} z_j/z_i;q)_{\theta_{i,N}}}.
\nonumber
\end{align}
This can be written explicitly as
\begin{align}\label{mrak}
&c_N(\theta;z_1,\ldots,z_N;q,\st)\\
\nonumber
&=
\prod_{k=2}^{N}
\prod_{1\leq i\leq j\leq k-1}
{(q^{\sum_{a=k+1}^N(\theta_{i,a}-\theta_{j+1,a})}\st z_{j+1}/z_i;q)_{\theta_{i,k}}\over
(q^{\sum_{a=k+1}^N(\theta_{i,a}-\theta_{j+1,a})}qz_{j+1}/z_i;q)_{\theta_{i,k}}}
{(q^{-\theta_{j,k}+\sum_{a=k+1}^N(\theta_{i,a}-\theta_{j,a})}qz_j/\st z_i;q)_{\theta_{i,k}}\over
(q^{-\theta_{j,k}+\sum_{a=k+1}^N(\theta_{i,a}-\theta_{j,a})}z_j/z_i;q)_{\theta_{i,k}}}=\\
\nonumber
&=\prod_{1\leq i<j\leq N}(q/\st)^{\theta_{i,j}}
\frac{(\st;q)_{\theta_{ij}}
(q^{\sum_{a=j+1}^N(\theta_{ia}-\theta_{ja})}\st z_j/z_i;q)_{\theta_{ij}}}
{(q;q)_{\theta_{ij}}
(q^{1+\sum_{a=j+1}^N(\theta_{ia}-\theta_{ja})}z_j/z_i;q)_{\theta_{ij}}}\times\\
\nonumber
& \prod_{k=3}^N\prod_{1\leq l<m<k}(q/\st)^{\theta_{l,k}}
\frac{(q^{\sum_{b=k+1}^N(\theta_{lb}-\theta_{mb})}\st z_m/z_l;q)_{\theta_{lk}}
(q^{-\theta_{lk}+\theta_{mk}-\sum_{b=k+1}^N(\theta_{lb}-\theta_{mb})}
\st z_l/z_m;q)_{\theta_{lk}}}
{(q^{1+\sum_{b=k+1}^N(\theta_{lb}-\theta_{mb})}z_m/z_l;q)_{\theta_{lk}}
(q^{1-\theta_{lk}+\theta_{mk}-\sum_{b=k+1}^N(\theta_{lb}-\theta_{mb})}
z_l/z_m;q)_{\theta_{lk}}}
\nonumber
\end{align}

%\begin{eg}
\ssec{eg}{Example}
\begin{align}
c_2&={(\st z_2/z_1;q)_{\theta_{1,2}} \over (qz_2/z_1;q)_{\theta_{1,2}} }
{(q^{-\theta_{1,2}}q/\st;q)_{\theta_{1,2}} \over (q^{-\theta_{1,2}};q)_{\theta_{1,2}} }
={(\st z_2/z_1;q)_{\theta_{1,2}} \over (qz_2/z_1;q)_{\theta_{1,2}} }
{(\st;q)_{\theta_{1,2}} \over (q;q)_{\theta_{1,2}} }(q/\st)^{\theta_{1,2}},\\
c_3&=
{(q^{\theta_{1,3}-\theta_{2,3}}\st z_2/z_1;q)_{\theta_{1,2}} \over
(q^{\theta_{1,3}-\theta_{2,3}}qz_2/z_1;q)_{\theta_{1,2}} }
{(q^{-\theta_{1,2}}q/\st;q)_{\theta_{1,2}} \over (q^{-\theta_{1,2}};q)_{\theta_{1,2}} }\\
&\times
{(\st z_2/z_1;q)_{\theta_{1,3}} \over (qz_2/z_1;q)_{\theta_{1,3}} }
{(q^{-\theta_{1,3}}q/\st;q)_{\theta_{1,3}} \over (q^{-\theta_{1,3}};q)_{\theta_{1,3}} }
{(\st z_3/z_1;q)_{\theta_{1,3}} \over (qz_3/z_1;q)_{\theta_{1,3}} }
{(q^{-\theta_{2,3}}qz_1/\st z_2;q)_{\theta_{1,3}} \over
 (q^{-\theta_{2,3}}z_1/z_2;q)_{\theta_{1,3}} }\nonumber\\
 &\times
 {(\st z_3/z_2;q)_{\theta_{2,3}} \over (qz_3/z_2;q)_{\theta_{2,3}} }
{(q^{-\theta_{2,3}}q/\st;q)_{\theta_{2,3}} \over (q^{-\theta_{2,3}};q)_{\theta_{2,3}} }.
\nonumber
\end{align}
%\end{eg}

\lem{sh4}
Let $\lambda=(\lambda_1,\lambda_2,\ldots)$ be a partition satisfying $\ell(\lambda)\leq N$.
The substitution $z_i=\st^{N-i}q^{\lambda_i}$ ($1\leq i\leq N$)
in $c_N(\theta;z_1,\ldots,z_N;q,\st)$ gives us the
coefficient $\psi_T$ in the tableau sum formula
\begin{eqnarray}
\psi_T(q,\st)=
c_N(\theta;\st^{N-1}q^{\lambda_1},\ldots,q^{\lambda_N};q,\st).
\end{eqnarray}
\elem

Let $y=(y_1,\ldots,y_N),z=(z_1,\ldots,z_N)$ be two sets of independent indeterminates.
Set
\begin{align}
z_i=\st^{N-i}q^{\lambda_i}\qquad (1\leq i\leq N).
\end{align}
For simplicity we use the notation
\begin{align}
y^{\lambda}=\prod_i y_i^{\lambda_i }.
\end{align}
Note that we have
\begin{align}
T_{q,y_i}y^{\lambda}= \st^{i-N}z_i  \cdot y^{\lambda}.
\end{align}

\defe{formal}
Define a formal power series $f_N(y,z;q,\st)\in y^{\lambda} \mathbb{F}(z)[[y_{i+1}/y_i,(i=1,\ldots,N-1)]]$ by
\begin{align}
f_N(y,z;q,\st)
=y^{\lambda} \sum_{\theta\in\mathsf{M}^{(N)}}c_N(\theta;z;q,\st)
\prod_{1\leq i<j \leq N} (y_j/y_i)^{\theta_{i,j}}. \label{blin}
\end{align}
\edefe

%%%%%%%%%%%%%%%%%%%%%%%%%%%%%%%%%%%%%%%%%%%%%%%%%
\ssec{tots}{Termination of the series $f_N(y,z;q,\st)$}
Let
$\lambda=(\lambda_1,\lambda_2,\ldots,\lambda_N)$
be a partition, while keeping $q,\st$ being generic.
Note that we have the following factor in the numerator of
$c_N(\theta;z_1,\ldots,z_N;q,\st)$:
\begin{align}
\prod_{k=1}^{N-1}
\prod_{1=1}^k
(q^{\sum_{a=k+1}^N(\theta_{i,a}-\theta_{i+1,a})}\st z_{i+1}/z_i;q)_{\theta_{i,k}}
=
\prod_{k=1}^{N-1}
\prod_{1=1}^k
(q^{\sum_{a=k+1}^N(\theta_{i,a}-\theta_{i+1,a})}q^{\lambda_{i+1}-\lambda_i};q)_{\theta_{i,k}}.
\end{align}
This vanishes unless the following set of inequalities are satisfied:
\begin{align}
0\leq \theta_{i,k}\leq \lambda_i-\lambda_{i+1}-\sum_{a=k+1}^N(\theta_{i,a}-\theta_{i+1,a})
\qquad (1\leq i<k\leq N).
\end{align}
Namely we have the vanishing of the coefficient $c_N(\theta;z_1,\ldots,z_N;q,\st)$'s
unless $\theta\in {\mathsf{Pol}}_\lambda \subset\mathsf{M}^{(N)}$.
Hence we find that
under the specialization in $z$, the infinite series $f_N(y,z;q,\st)$ terminates into a finite one.

\prop{tokusyuka}
Let
$\lambda=(\lambda_1,\lambda_2,\ldots,\lambda_N)$
be a partition, and set $z_i=\st^{N-i}q^{\lambda_i}$. Then we have
\begin{align}
f_N(y,z;q,\st)
&=y^\lambda \sum_{\theta\in{\mathsf{Pol}}_\lambda}c_N(\theta;z;q,\st)
\prod_{1\leq i<j \leq N} (y_j/y_i)^{\theta_{i,j}}\\
&=\sum_{T} \psi_T(q,\st)y^T=P_\lambda(y,q,\st). \nonumber
\end{align}
\eprop

\prop{dai-ichi}
Let $y=(y_1,\ldots,y_N)$ and $z=(z_1,\ldots,z_N)$  be generic.
We have
\begin{eqnarray}
\sD^1_{N,y} f_N(y,z;q,\st)=\sum_{i=1}^N z_i \cdot f_N(y,z;q,\st).\label{dai-ichi-shiki}
\end{eqnarray}
\eprop

\lem{takoushiki}
Let
$u(z_1,\ldots,z_N)\in \mathbb{F}[z_1,z_2,\ldots,z_N]$.
If we have $u(\st^{N-1}q^{\lambda_1},\st^{N-2}q^{\lambda_2},\ldots,q^{\lambda_N})=0$
for any partition $\lambda=(\lambda_1,\lambda_2,\ldots,\lambda_N)$, then $u(z_1,\ldots,z_N)=0$.
\elem

\prf
We prove this by the induction on $N$.
When $N=1$, it is true.
Assume it holds for $N-1$.
Expand
$u(z_1,\ldots,z_N)=\sum_k u_k(z_2,\ldots,z_N)z_1^k$.
Fix $\lambda_2,\ldots,\lambda_N$ and vary $\lambda_1(\geq\lambda_2)$,
then all the coefficients of $z_1^k$, i.e.
$u_k(\st^{N-2}q^{\lambda_2},\st^{N-3}q^{\lambda_3},\ldots,q^{\lambda_N})$
should vanish. Now we let  $\lambda_2,\ldots,\lambda_N$ vary and
conclude that $u_k(z_2,\ldots,z_N)=0$ by the assumption.
\epr

\noindent
{\em Proof of~\refp{dai-ichi}}.
Set
\begin{align*}
&y^{-\lambda}(\mbox{LHS(\ref{dai-ichi-shiki})}-\mbox{RHS(\ref{dai-ichi-shiki})})\\
&=\sum_{k_1,\ldots,k_{N-1}\geq0} r_{k_1,\ldots,k_{N-1}}(z) \prod_{i=1}^{N-1}(y_{i+1}/y_i)^{k_i}
\in  \mathbb{F}(z)[[y_{i+1}/y_i,(i=1,2,\ldots,N-1)]].
\end{align*}
{}From~\refp{tokusyuka} and~\refl{takoushiki},
we have
 $r_{k_1,\ldots,k_{N-1}}(z) =0$ for all $k_1,\ldots,k_{N-1}\geq0$.
\qed

\sec{van}{Vanishing}

\ssec{som}{Sommese vanishing}
We need the following version of Sommese vanishing theorem.
Let $p:\ X\to Y$ be a flat morphism between smooth projective complex
varieties. Let $\CL$ be a line bundle on $X$ whose restriction to every
fiber of $p$ is $l$-ample~\cite[Definition~6.5]{EV} for certain $l\in\BN$.
Also, suppose the Iitaka dimension $\kappa(\CL_y)$~\cite[Definition~5.3]{EV}
of the restriction
of $\CL$ to every fiber $X_y=p^{-1}(y),\ y\in Y$, equals $\dim X_y=\dim X-
\dim Y$. Finally, suppose $\dim X-\dim Y-l>M$ for some $M\in\BN$.

\th{somm}
Under the above assumptions, $H^i(X,\Omega^j_X\otimes\CL^{-1})=0$ for
$i+j<M$.
\eth

\prf
By the Leray spectral sequence, it suffices to prove
$R^ip_*(\Omega^j_X\otimes\CL^{-1})=0$ for $i+j<M$. First, we restrict to
a nonempty open $U\subset Y$ over which $p$ is smooth. We set $X_U=p^{-1}(U)$.
Then $\Omega^j_{X_U}$ has a filtration whose associated graded bundle is a
direct sum of the sheaves $\Omega^k_{X_U/U}\otimes p^*\Omega^{j-k}_U$ over
$k\leq j$. Here $\Omega^k_{X_U/U}$ is the bundle of relative $k$-forms.
By the projection formula it suffices to prove
$R^ip_*(\Omega^k_{X_U/U}\otimes\CL^{-1})=0$ for $i+k<M$.
By the base change, it suffices to know for any $y\in U$ that
$H^i(X_y,\Omega^k_{X_y}\otimes\CL^{-1})=0$ for $i+k<M$. But this is nothing
but Sommese vanishing~\cite[Corollary~6.6]{EV} on the smooth projective
variety $X_y$. So we conclude
$R^ip_*(\Omega^j_X\otimes\CL^{-1})|_U=0$ for $i+j<M$.

Now to prove $R^ip_*(\Omega^j_X\otimes\CL^{-1})=0$ for $i+j<M$ it suffices
to know that $R^ip_*(\Omega^j_X\otimes\CL^{-1})$ has no torsion for $i+j<M$.
We will prove this by induction in $\dim Y$ and
the dimension of the support of torsion.
Let $Z\subset Y$ be the support of $R^ip_*(\Omega^j_X\otimes\CL^{-1})$.
Suppose $\dim Z>0$. Then according to Kleiman's generic transversality
theorem~\cite{K} there exists a hyperplane section $Y'\subset Y$ intersecting
$Z$ transversally at a smooth point $z\in Z$ and such that $X':=p^{-1}(Y')$
is smooth. By the base change, the support of
$R^ip_*(\Omega^j_X\otimes\CL^{-1}|_{X'})$ contains $Z'$ defined as the
irreducible component of $Z\cap Y'$ containing $z$.
However, we have an exact sequence of vector bundles on $X'$:
$$0\to \CN^*_{X'/X}\otimes\Omega^{j-1}_{X'}\to\Omega^j_X|_{X'}\to\Omega^j_{X'}\to0$$
and the conormal bundle $\CN^*_{X'/X}=p^*\CN^*_{Y'/Y}$. By the projection formula
and by the induction (in $\dim Y$) assumption we have
$R^ip_*(\CN^*_{X'/X}\otimes\Omega^{j-1}_{X'}\otimes\CL^{-1})=0=
R^ip_*(\Omega^j_{X'}\otimes\CL^{-1})$ for $i+j<M$. Hence
$R^ip_*(\Omega^j_X|_{X'}\otimes\CL^{-1})=0$ which contradicts to
$Z'\ne\emptyset$.

It remains to establish the base of induction: $\dim Z=0$.
We choose the minimal $i_0$ among all $i$ such that
$R^ip_*(\Omega^j_X\otimes\CL^{-1})\ne0$ for some $j$ such that $i+j<M$.
Let $y\in Z\subset Y$ be a point in the (finite) support of
$R^{i_0}p_*(\Omega^j_X\otimes\CL^{-1})$. Let us choose a sufficiently ample
line bundle $\CM$ on $Y$. Then by the projection formula
$R^ip_*(\Omega^j_X\otimes(\CL\otimes p^*\CM)^{-1})=
R^ip_*(\Omega^j_X\otimes\CL^{-1})\otimes\CM^{-1}$, and by the Leray spectral
sequence $H^{i_0}(X,\Omega^j_X\otimes(\CL\otimes p^*\CM)^{-1})\ne0$ (the LHS
contains a direct summand
$R^{i_0}p_*(\Omega^j_X\otimes\CL^{-1})_y\otimes\CM^{-1}_y$). However, by the
Sommese vanishing~\cite[Corollary~6.6]{EV} applied to the line bundle
$\CL\otimes p^*\CM$ on $X$ (with $\CM$ sufficiently ample), we must have
$H^{i_0}(X,\Omega^j_X\otimes(\CL\otimes p^*\CM)^{-1})=0$. This contradiction
proves we cannot have $\dim Z=0$.

This completes the proof of the theorem.
\epr

\ssec{parab}{Parabolic Laumon spaces}
Recall the notations of~\cite{BF12}. We denote by $\QM^\alpha$ the
Drinfeld moduli space of degree $\alpha$ quasimaps from $\bC\simeq\BP^1$
to the flag variety $\CB=G/B$ of $G=\SL(N)$. Here
$\alpha=(d_1,\ldots,d_{N-1})\in\BN^{N-1}$. We denote by
$\pi_\alpha:\ \CQ^\alpha\to\QM^\alpha$ the Laumon resolution of
$\QM^\alpha$~\cite{Ku}. Given a subminimal parabolic (with Levi of semisimple
rank 1)
$\SL(N)\supset P=P_i\supset B$ we consider the corresponding parabolic
Laumon space $\CQ^\balpha_P$ (see e.g.~\cite{BF10}), and the natural
projection $\varpi_\alpha:\ \CQ^\alpha\to\CQ^\balpha_P$. Here
$\balpha:=(d_1,\ldots,d_{i-1},d_{i+1},\ldots,d_{N-1})$.

Recall~\cite{BF12} that $V_{\check\omega_i}=\Lambda^i\BC^N,\ 1\leq i\leq N-1$,
are the fundamental $\SL(N)$-modules,
and $\QM^\alpha$ is equipped with a closed embedding
$\psi_\alpha:\ \QM^\alpha\hookrightarrow
\prod_{i\in I}\BP\Gamma(\bC,V_{\check\omega_i}\otimes\CO(\langle\alpha,
\check\omega_i\rangle))$.
Given an $\SL(N)$-weight $\check\lambda=
\sum_{i\in I}d_i\check\omega_i\in\Lambda^\vee$ we define a line bundle
$\CO(\check\lambda)^\alpha$ on $\QM^\alpha_\fg$ as
$\psi_\alpha^*\bigotimes_{i\in I}\CO(d_i)$.
Suppose $\chlam$ is {\em not}
dominant, i.e. $l_i<0$ for some $1\leq i\leq N-1$. We fix such an $i$ from
now on, and we set $\CL:=\pi_\alpha^*\CO(-\chlam)$. For $y\in\CQ^\balpha_P$
we denote by $X_y$ the fiber $\varpi_\alpha^{-1}(y)$ with the reduced scheme
structure. Our aim is to study the ampleness properties of the line bundle
$\CL_y:=\CL|_{X_y}$. They are summarized in the following

\prop{parabol}
(a) $\CL_y$ is generated by the global sections, and gives rise to a
morphism $\phi:\ X_y\to\BP(\Gamma^*(X_y,\CL_y))$. We denote by
$\overline{X}_y$ the image of $\phi$ (with the reduced closed subscheme
structure).

(b) The morphism $X_y\stackrel{\phi}{\longrightarrow}\overline{X}_y$
equals $X_y\stackrel{\pi_\alpha}{\longrightarrow}\pi_\alpha(X_y)$, where
$\pi_\alpha(X_y)\subset\QM^\alpha$ is equipped with the reduced closed
subscheme structure.

(c) For a fixed $\balpha$ and $d_i\gg0$ we have
$\dim X_y=\dim\overline{X}_y=2d_i-d_{i-1}-d_{i+1}+1$; in particular,
$\varpi_\alpha$ is flat.

(d) Let $l_y:=\max\{\dim\phi^{-1}(z), z\in\overline{X}_y\}$. For a fixed
$\balpha$ and $M\in\BN$, there exists $D_i$ such that for
$d_i>D_i$ and any $y\in\CQ^\balpha_P$ we have $\dim X_y-l_y>M$.
\eprop

\prf
(a) and (b) are clear from definitions.
A point $y\in\CQ^\balpha_P$ is represented by a collection of
locally free subsheaves $0=\CW_0\subset\CW_1\subset\ldots\subset\CW_{i-1}\subset
\CW_{i+1}\subset\ldots\subset\CW_{N-1}\subset\CW_N=\CO^N_\bC$ such that
$\operatorname{rk}\CW_j=j$, and $\deg\CW_j=-d_j$.
The fiber $X_y$ is the moduli space of subsheaves
$\overline\CW_i\subset\CW^{i+1}_{i-1}:=\CW_{i+1}/\CW_{i-1}$ of generic rank 1
and degree $d_{i-1}-d_i$. For such a sheaf $\overline\CW_i$ we denote
by $\overline{\underline\CW}{}_i$ its {\em saturation} i.e. the maximal
subsheaf of $\CW^{i+1}_{i-1}$ containing $\overline\CW_i$, of generic rank 1,
and such that $\CW^{i+1}_{i-1}/\overline{\underline\CW}{}_i$ has no torsion.
We also define the {\em defect} $\on{def}\overline\CW_i$ as the cycle
of the torsion sheaf $\overline{\underline\CW}{}_i/\overline\CW_i$.
Two points $\overline\CW_i,\overline\CW{}'_i$ are in the same fiber of
$\phi=\pi_\alpha|_{X_y}$ iff their saturations and defects coincide.
In particular,
$\phi$ is one-to-one when restricted to the open subset $U\subset X_y$
formed by all the saturated $\overline\CW_i$.

To prove (c) we must check
that $U$ is nonempty for $d_i\gg0$. This is evident. To finish the proof of (c)
it remains to compute $\dim U$. Let us decompose $\CW^{i+1}_{i-1}\simeq
(\CW^{i+1}_{i-1})^{\on{tors}}\oplus(\CW^{i+1}_{i-1})^{\on{free}}$ into a direct sum
of a torsion sheaf and a locally free sheaf. Then a point of $U$ is
represented by $\overline\CW_i\simeq(\CW^{i+1}_{i-1})^{\on{tors}}\oplus
\overline\CW{}_i^{\on{free}}$ where $\overline\CW{}_i^{\on{free}}\subset
(\CW^{i+1}_{i-1})^{\on{free}}$ is a line subbundle of degree
$d_{i-1}-d_i-\dim(\CW^{i+1}_{i-1})^{\on{tors}}$. Locally around $\overline\CW_i$,
$U$ is isomorphic to
$\BP\Hom(\overline\CW{}_i^{\on{free}},\CW^{i+1}_{i-1})$.
For $d_i\gg0$ the latter space is $\BP^{2d_i-d_{i-1}-d_{i+1}+1}$ which
completes the proof of (c).

To prove (d) we fix a saturated subsheaf $\overline{\underline\CW}{}_i=
\overline{\underline\CW}{}_i^{\on{tors}}\oplus
\overline{\underline\CW}{}_i^{\on{free}}$,
and a defect $\delta\in\bC^{(d)}$. We have to estimate the dimension of the
moduli space $\phi^{-1}(z)$ of subsheaves
$\overline\CW_i\subset\CW^{i+1}_{i-1}$ with given
saturation $\overline{\underline\CW}{}_i$ and $\on{def}\overline\CW_i=\delta$.
If $\overline\CW_i^{\on{prfr}}$ stands for the image of projection of
$\overline\CW_i$ to $\overline{\underline\CW}{}_i^{\on{free}}$ along
$\overline{\underline\CW}{}_i^{\on{tors}}$, then there are finitely many
possible values of
$\overline\CW_i^{\on{prfr}}\subset\overline{\underline\CW}{}_i^{\on{free}}$;
more precisely, not more than $\tau^\tau$ where
$\tau=\dim(\CW^{i+1}_{i-1})^{\on{tors}}$ (the only ambiguity in the choice
of $\overline\CW_i^{\on{prfr}}\subset\overline{\underline\CW}{}_i^{\on{free}}$
can occur at the support of $(\CW^{i+1}_{i-1})^{\on{tors}}$).
Now for a fixed value of
$\overline\CW_i^{\on{prfr}}\subset\overline{\underline\CW}{}_i^{\on{free}}$
the dimension of the corresponding stratum of the moduli space in question
is independent of $d_i$. Hence, with $d_i$ growing, $\dim X_y-\dim\phi^{-1}(z)$
grows uniformly in $y$ and $z$.

The proposition is proved.
\epr

\ssec{vani}{Vanishing Theorem}
We combine~\refp{parabol} and~\reft{somm} setting
$Y=\CQ^\balpha_P,\ X=\CQ^\alpha,\ p=\varpi_\alpha,\
\CL=\CO(-\chlam)$. We arrive at the following

\th{vanis}
Let $\chlam=\sum_{i=1}^{N-1}l_i\check\omega_i$ be a non-dominant weight, i.e.
$l_i<0$ for certain $1\leq i\leq N-1$.
We fix $j,k\in\BN$ and $\balpha=(d_1,\ldots,d_{i-1},d_{i+1},\ldots,d_{N-1})$.
Then for $d_i\gg0$ we have
$H^k(\CQ^\alpha,\Omega^j_{\CQ^\alpha}\otimes\CO(\chlam))=0. \qed$
\eth

\sec{DR}{Euler characteristics of twisted De Rham complexes}

\ssec{gen}{Generating functions}
For a weight $\chlam$ we consider
$\chi(H^\bullet(\CQ^\alpha,\Omega^\bullet_{\CQ^\alpha}\otimes\CO(\chlam)))$
as a virtual graded $\BG_m\times T$-module.
Here $T$ is the Cartan torus of $\SL(N)$,
and the grading is via De Rham degree $\Omega^\bullet_{\CQ^\alpha}$.
The generating function of its character is
$[H^\bullet(\CQ^\alpha,\Omega^\bullet_{\CQ^\alpha}\otimes\CO(\chlam))]:=
\sum_{i,j}(-1)^{i+j}t^j[H^i(\CQ^\alpha,\Omega^j_{\CQ^\alpha}\otimes\CO(\chlam))]$
a function of $q,t,z$ where $q$ is the coordinate on $\BG_m$, and
$z$ are the coordinates on $T$.
We define $H_\chlam(q,t,z)$ as the limit of
$[H^\bullet(\CQ^\alpha,\Omega^\bullet_{\CQ^\alpha}\otimes\CO(\chlam))]$ as
$\alpha\to\infty$. We will see that the limit exists as a formal series in
$q,t,z$ converging to a rational function in $q,t,z$. For instance, if
$\chlam$ is {\em not} dominant, then according to~\reft{vanis},
$H_\chlam(q,t,z)=0$.

\ssec{coh}{Betti cohomology of Laumon spaces}
We start with a computation of $H_0(q,t,z)$.

\prop{h0}
$H_0(q,t,z)=\frac{(1+t)(1+t+t^2)\ldots(1+t+\ldots+t^{N-1})}
{(1-t^{N-1})^2(1-t^{N-2})^4\ldots(1-t^3)^{2N-6}(1-t^2)^{2N-4}}\cdot
\frac{1}{(1-t^N)(1-t)^{N-2}}$
\eprop

\prf
According to~\cite[Theorem~2.9]{FK}, the Betti cohomology
$H^\bullet(\CQ^\alpha,\BC)$ carries a Tate Hodge structure, so
$H^i(\CQ^\alpha,\Omega^j_{\CQ^\alpha})=0$ unless $i=j$, while
$H^i(\CQ^\alpha,\Omega^i_{\CQ^\alpha})=H^{2i}(\CQ^\alpha,\BC)$.
The action of $\BG_m\times T$ on the latter space is clearly trivial,
so $H_0(q,t,z)=H_0(t)$ is the $\alpha\to\infty$ limit of Poincar\'e
polynomials $P_\alpha(t):=\sum_it^i\dim H^{2i}(\CQ^\alpha,\BC)$.
Now $P_\alpha(t)$ is calculated in~\cite[Theorem~2.7]{FK} (under the perverse
normalization). The $\alpha\to\infty$ limit $P_\infty(t)$ is the product
$W(t)\cdot F(t)$ where $W(t)$ is the Poincar\'e polynomial of the flag
variety $\CB$, that is $W(t)=N!_t=(1+t)(1+t+t^2)\ldots(1+t+\ldots+t^{N-1})$.
Furthermore, $F(t)=\sum F_it^i$ where $F_i$ is the number of unordered
collections of positive roots $\alpha_1,\ldots,\alpha_k,\beta_1,\ldots,\beta_m
\in R^+(\mathfrak{sl}_N)$ such that none of $\alpha_1,\ldots,\alpha_k$ is
simple, and $\sum_{j=1}^k(|\alpha_j|-1)+\sum_{l=1}^m(|\beta_l|+1)=i$.
Here $|\beta|:=(\beta,\rho)$.

The proposition follows.
\epr

\ssec{loc}{Local Laumon spaces}
Recall that $\fQ^\alpha\subset\CQ^\alpha$ is a locally closed {\em local}
Laumon moduli space of quasiflags based at $\infty\in\bC$,
see e.g.~\cite{FFFR}. Similarly to~\refss{gen} we introduce the generating
function $\fJ_\alpha(q,t,z)=[H^\bullet(\fQ^\alpha,\Omega^\bullet_{\fQ^\alpha})]:=
\sum_{i,j}(-1)^{i+j}t^j[H^i(\fQ^\alpha,\Omega^j_{\fQ^\alpha})]$.
We compute the Euler characteristic of
$H^\bullet(\fQ^\alpha,\Omega^\bullet_{\fQ^\alpha})$ via the Atiyah-Bott-Lefschetz
localization to the fixed points of $\BG_m\times T$ in $\fQ^\alpha$.
The characters of $\BG_m\times T$ in the tangent spaces of the fixed points
are computed in~\cite[Proposition~2.18a]{BF0}. To write down the answer we
recall the necessary notation. The fixed points are numbered by the collections
$\widetilde{d}=(d_{ij})_{N-1\geq i\geq j\geq1}$ such that $d_{ij}\leq d_{kj}$
for $i\geq k\geq j$, and $d_{i,1}+d_{i,2}+\ldots+d_{i,i}=d_i$ (recall that
$\alpha=(d_1,\ldots,d_{N-1})$). For $1\leq k<l\leq N$ we set
$\theta_{kl}=d_{l-1,k}-d_{lk}$ where $d_{Nk}:=0$. Conversely,
$d_{ij}=\theta_{j,i+1}+\theta_{j,i+2}+\ldots+\theta_{j,N}$.
%Recall that the set $\{(k,l):\ 1\leq k<l\leq N\}$ is denoted by $\sM^{(N)}$.
Recall the set $\sM^{(N)}$ introduced in~(\ref{215}).
Let $x_i$ stand for the character of the dual torus
$\check T$ corresponding to the simple coroot $\alpha_i$. For
$\alpha\in\BN^{N-1}$ the corresponding character of $\check T$ is denoted by
$x^\alpha$. We consider the formal generating function
$$
J(q,t,z,x)=\sum_{\alpha\in\BN^{N-1}}x^\alpha\fJ_\alpha(q,t,z).
$$
Now
\eq{60}
J(q,t,z,x)=\sum_{(\theta_{ij})\in\sM^{(N)}}C_{(\theta_{ij})}(q,t,z)
\prod_{1\leq i<j\leq N}
(x_i\cdots x_{j-1})^{\theta_{ij}}
\end{equation}
where
\begin{align}
\label{shira}
& C_{(\theta_{ij})}(q,t,z)=\prod_{1\leq i<j\leq N}
\frac{(qt;q)_{\theta_{ij}}
(q^{1+\sum_{a=j+1}^N(\theta_{ia}-\theta_{ja})}tz_j/z_i;q)_{\theta_{ij}}}
{(q;q)_{\theta_{ij}}
(q^{1+\sum_{a=j+1}^N(\theta_{ia}-\theta_{ja})}z_j/z_i;q)_{\theta_{ij}}}\times\\
& \prod_{k=3}^N\prod_{1\leq l<m<k}
\frac{(q^{1+\sum_{b=k+1}^N(\theta_{lb}-\theta_{mb})}tz_m/z_l;q)_{\theta_{lk}}
(q^{1-\theta_{lk}+\theta_{mk}-\sum_{b=k+1}^N(\theta_{lb}-\theta_{mb})}
tz_l/z_m;q)_{\theta_{lk}}}
{(q^{1+\sum_{b=k+1}^N(\theta_{lb}-\theta_{mb})}z_m/z_l;q)_{\theta_{lk}}
(q^{1-\theta_{lk}+\theta_{mk}-\sum_{b=k+1}^N(\theta_{lb}-\theta_{mb})}
z_l/z_m;q)_{\theta_{lk}}}
\nonumber
\end{align}
--- this is a restatement of~\cite[Proposition~2.18a]{BF0}.

\ssec{loc stab}{Local stabilization}
We will prove that as $\alpha\to\infty$,
the series $\fJ_\alpha(q,t,z)$ tends to the limit $\fJ_\infty(q,t,z)$.
More precisely, for any $n,m\in\BN$ the coefficient of $q^nt^m$ in
$\fJ_\alpha(q,t,z)$ for $\alpha\gg0$ is independent of $\alpha$. The resulting
series will be denoted by $\fJ_\infty(q,t,z)$. The existence of the limit and
computation of its value follows from~\reft{shir} below
and~\cite[Proposition~3.11]{CE}. For the reader's convenience we present a
more elementary computation of the limit.

We introduce a function $z_N$ on the Cartan torus $T$ of $\SL(N)$ defined
as $z_N:=z_1^{-1}\cdots z_{N-1}^{-1}$.

\th{junichi}
$$\lim_{\alpha\to\infty}\fJ_\alpha(q,t,z)=\prod_{1\leq i<j\leq N}
{(qt z_j/z_i;q)_{\infty } \over (q z_j/z_i;q)_{\infty}}\times
\left({(qt;q)_{\infty } \over (q;q)_{\infty}}\right)^{N-1}\times
\prod_{i=1}^{N-2}
\left({(qt^{i+1};q)_\infty \over (t^i;q)_\infty}\right)^{N-i-1}.$$
\eth

\prf
Let $\alpha=\sum_{i=1}^{N-1}\ell_i\alpha_i, \ell_i\in\BN$.
We study the stabilization of $\fJ_\alpha(q,t,z)$ in the sector
$\ell_1\gg \ell_2 \gg \cdots\gg \ell_{N-1}\gg 0$.
Note that we can recast the coefficient $C_{(\theta_{i,j})} (q,t,z)$ as follows:
\begin{align*}
C_{(\theta_{i,j})} (q,t,z)&=
\prod_{1\leq i<j\leq N}
{(qt;q)_{\theta_{i,j}} \over (q;q)_{\theta_{i,j}}}
{(qt z_j/z_i;q)_{\theta_{i,N}} \over (q z_j/z_i;q)_{\theta_{i,N}}}\times
\prod_{k=3}^{N} F_{k},\\
 F_{k}&=
 \prod_{1\leq l<m< k}
 {(q^{\sigma_{l,k}-\sigma_{m,k}} qt z_m/z_l;q)_{\theta_{l,k-1}} \over (q^{\sigma_{l,k}-\sigma_{m,k} }q z_m/z_l;q)_{\theta_{l,k-1}} }
  {(q^{-\sigma_{l,k}+\sigma_{m,k}}q t z_l/z_m;q)_{\theta_{l,k}} \over (q^{-\sigma_{l,k}+\sigma_{m,k} }q z_l/z_m;q)_{\theta_{l,k}} },\\
  \sigma_{l,k}&=\sum_{b=k}^N \theta_{l,b}.
\end{align*}

Define $F^{(n)}$ ($0\leq n\leq N-1$) recursively by setting $F^{(0)}=C_{(\theta_{i,j})} (q,t,z)$,
\begin{align*}
F^{(1)}&=\lim_{\ell_{1}\rightarrow \infty}
\sum_{\theta_{1,2}\geq 0 \atop \ell_1=\sigma_{1,2}}F^{(0)},\qquad
F^{(2)}= \lim_{\ell_{2}\rightarrow \infty}
\sum_{\theta_{1,3},\theta_{2,3}\geq 0 \atop \ell_2=\sigma_{1,3}+\sigma_{2,3}} F^{(1)},\ldots,
\\
F^{(n)}&=\lim_{\ell_{n}\rightarrow \infty}
\sum_{\theta_{1,n+1},\ldots,\theta_{n,n+1}\geq 0 \atop \ell_{n}=\sigma_{1,n+1}+\cdots+\sigma_{n,n+1}}F^{(n-1)},\ldots,\\
F^{(N-1)}&=\lim_{\ell_{N-1}\rightarrow \infty}
\sum_{\theta_{1,N},\ldots,\theta_{N-1,N}\geq 0 \atop \ell_{N-1}=\theta_{1,N}+\cdots+\theta_{N-1,N}}F^{(N-2)}.
\end{align*}
Then the coefficient we are interested in is $F^{(N-1)}$.
The stabilization of   $F^{(N-1)}$ can be studied and stated explicitly as follows.

\lem{prp}
 For $n=1,\ldots,N-2$, we have
 \begin{align*}
 F^{(n)}&=
 \prod_{1\leq i<j\leq n+1}
{(qt;q)_{\infty } \over (q;q)_{\infty}}
{(qt z_j/z_i;q)_{\theta_{i,N}} \over (q z_j/z_i;q)_{\theta_{i,N}}}\times \\
&\times
\left({(q;q)_\infty \over (qt;q)_\infty }\right)^{\left(n\atop 2\right)} \prod_{i=1}^{n-1}
\left({(qt^{i+1};q)_\infty \over (t^i;q)_\infty}\right)^{n-i}\times\\
&\times
\sum_{(\theta_{i,j} )_{j\geq n+2} \atop
\ell_i=\sigma_{1,i+1}+\cdots+\sigma_{i,i+1}~(n+1\leq i\leq N-1)}
\prod_{1\leq i<j\leq N \atop n+2\leq j}
{(qt;q)_{\theta_{i,j}} \over (q;q)_{\theta_{i,j}}}
{(qt z_j/z_i;q)_{\theta_{i,N}} \over (q z_j/z_i;q)_{\theta_{i,N}}}\times
\prod_{k=n+2}^{N}  F_{k},\\
 \end{align*}
  and
  \begin{align*}
  F^{(N-1)}=
 \prod_{1\leq i<j\leq N}
{(qt z_j/z_i;q)_{\infty } \over (q z_j/z_i;q)_{\infty}}\times
\left({(qt;q)_{\infty } \over (q;q)_{\infty}}\right)^{N-1}\times
\prod_{i=1}^{N-2}
\left({(qt^{i+1};q)_\infty \over (t^i;q)_\infty}\right)^{N-i-1}.
 \end{align*}
\elem

\prf
 We prove the statement by the recursive use of the summation formula
associated with the
root lattice of type $A_n$~\cite[the table at p.~136 and references therein]{I}:
\begin{align*}
\sum_{\chi\in Q} \prod_{\alpha\in R} {(q^{1+\langle \alpha,\chi \rangle }t z_\alpha;q)_\infty \over
(q^{1+\langle \alpha,\chi \rangle }z_\alpha;q)_\infty}
&=
\prod_{\alpha>0}
{(q t^{\langle \rho,\alpha\rangle+1};q)_\infty (q^{\delta_\alpha}t^{\langle \rho,\alpha\rangle-1};q)_\infty
\over
(q t^{\langle \rho,\alpha\rangle};q)_\infty
(t^{\langle \rho,\alpha\rangle};q)_\infty }\\
&=
\left({(q;q)_\infty \over (qt;q)_\infty }\right)^{n-1} \prod_{i=1}^{n-1}{(qt^{i+1};q)_\infty \over (t^i;q)_\infty},
\end{align*}
where $2 \rho=\sum_{\alpha>0}\alpha$,
$\delta_\alpha=1$ if $\alpha$ is a simple root and $\delta_\alpha=0$ otherwise.

It is clear that we have $F^{(1)}$ by taking the limit $\ell_1\rightarrow \infty$, namely  letting $\theta_{1,2}\rightarrow \infty$
while fixing all the other $\theta_{i,j}$'s.
The passage from $F^{(n)}$ to $F^{(n+1)}$ can be studied as follows.
We need to take the limit $\ell_n\rightarrow \infty$
and perform the $(n-1)$-dimensional summation with respect to the
variables $\theta_{1,n+1},\ldots,\theta_{n,n+1}$ with the constraint
$\ell_n=\sigma_{1,n+1}+\cdots+\sigma_{n,n+1}$.
It can be easily shown that,
the most dominating terms, as a Taylor series in $q$ and $t$, come from the vicinity of $\sigma_{l,n+1}-\sigma_{m,n+1}\sim 0$ ($1\leq l< m\leq n $).
Therefore the dominating contributions come from $\theta_{1,n+1},\ldots,\theta_{n,n+1}\gg 0 $.
Hence the $(n-1)$-dimensional summation
can be taken by using the above mentioned summation formula for type $A_{n}$.
\epr

This completes the proof of the theorem.
\epr

\ssec{loc to glob}{From local to global Laumon spaces}
The following lemma is very similar
to~\cite[Lemma~4.2]{BF12}:

\lem{bf12}
$$[H^\bullet(\CQ^\alpha,\Omega^\bullet_{\CQ^\alpha}\otimes\CO(\chlam))]=
\sum_{\substack{\gamma+\beta=\alpha\\ w\in W}}
z^{w\check\lambda}q^{\langle\gamma,\check\lambda\rangle}
\fJ_\gamma(q^{-1},t,wz)\fJ_\beta(q,t,wz)
\prod_{\check\alpha\in\check{R}{}^+}\frac{1-twz^{\check\alpha}}
{1-wz^{\check\alpha}}.$$
\elem

\prf
Atiyah-Bott-Lefschetz localization to the fixed points of $\BG_m\times T$
in $\CQ^\alpha$, see~\cite[Proof of~Theorem~5.8]{FFFR}.
\epr

\ssec{glob stab}{Global stabilization}
We consider the formal generating function
$\fJ(q,t,z,x)=\prod_{i=1}^{N-1}x_i^{\log(\check\omega_i)/\log q}J(q,t,z,x)$.

Note that if we plug $x=q^{\check\lambda}$ into $J(q^{-1},t,z,x)$ or into
$\fJ(q^{-1},t,z,x)$, then for a dominant weight $\check\lambda$ these formal
series converge, and we have $\fJ(q^{-1},t,z,q^{\check\lambda}):=\prod_{i=1}^{N-1}
(q^{\langle\alpha_i,\check\lambda\rangle})^{\log(\check\omega_i)/\log q}
J(q^{-1},t,z,q^{\check\lambda})=
z^{\check\lambda}J(q^{-1},t,z,q^{\check\lambda})$ (a formal Taylor series in
$q$ and $t$ with coefficients in Laurent polynomials in $z$).

Recall the generating function $H_\chlam(q,t,z)$ introduced in~\refss{gen}.
The following proposition is very similar to~\cite[Proposition~4.4]{BF12}:

\prop{brav}
$$
H_\chlam(q,t,z)=
\sum_{w\in W}\fJ(q^{-1},t,wz,q^{\check\lambda})\fJ_\infty(q,t,wz)
\prod_{\check\alpha\in\check{R}{}^+}\frac{1-twz^{\check\alpha}}
{1-wz^{\check\alpha}}.
$$
\eprop

\prf
As $\alpha$ goes to $\infty$, the formula of~\refl{bf12} goes to
$$\sum_{\substack{\gamma\in\Lambda_+\\ w\in W}}
z^{w\check\lambda}q^{\langle\gamma,\check\lambda\rangle}
\fJ_\gamma(q^{-1},t,wz)\fJ_\infty(q,t,wz)
\prod_{\check\alpha\in\check{R}{}^+}\frac{1-twz^{\check\alpha}}
{1-wz^{\check\alpha}}=$$
\eq{middle}
\sum_{w\in W}
z^{w\check\lambda}J(q^{-1},t,wz,q^{\check\lambda})\fJ_\infty(q,t,wz)
\prod_{\check\alpha\in\check{R}{}^+}\frac{1-twz^{\check\alpha}}
{1-wz^{\check\alpha}}=
\end{equation}
$$\sum_{w\in W}\fJ(q^{-1},t,wz,q^{\check\lambda})\fJ_\infty(q,t,wz)
\prod_{\check\alpha\in\check{R}{}^+}\frac{1-twz^{\check\alpha}}
{1-wz^{\check\alpha}}.$$
\epr

\sec{DE}{Difference equations}

\ssec{eu}{Euler characteristics of De Rham complexes of local Laumon spaces}
For $1\leq i\leq N$, we consider the difference operator $T_{i,q^{\pm1}}$
on functions of $q,t,z,x$
defined as follows: $T_{i,q^{\pm1}}F(q,t,z,x_1,\ldots,x_{N-1}):=
F(q,t,z,x_1,\ldots,x_{i-2},q^{\mp1}x_{i-1},q^{\pm1}x_i,x_{i+1},\ldots,x_{N-1})$.
We define
$$\sD:=\sum_{i=1}^N\prod_{j<i}
\frac{1-q^{-1}t^{i-j-1}x_j\cdots x_{i-1}}{1-t^{i-j}x_j\cdots x_{i-1}}
\prod_{k>i}\frac{1-qt^{k-i+1}x_i\cdots x_{k-1}}{1-t^{k-i}x_i\cdots x_{k-1}}
T_{i,q^{-1}}$$
Recall the function $z_N$ on the Cartan torus $T$ of $\SL(N)$ defined
as $z_N:=z_1^{-1}\cdots z_{N-1}^{-1}$.

\th{shir}
$\sD\fJ(q,t,z,x)=(z_1+\ldots+z_N)\fJ(q,t,z,x)$.
\eth

\prf
We just recall~\refp{dai-ichi} and compare~\refe{60} and~(\ref{shira}) with
formula~(\ref{mrak}) for an
eigenfunction $f_N(q,\st,z_1,\ldots,z_N,y_1,\ldots,y_N)$ of the
difference operator
$$\sD^1_N=\sum_{i=1}^Nz_i\prod_{j<i}\frac{1-\st^{-1}y_i/y_j}{1-y_i/y_j}
\prod_{k>i}\frac{1-\st y_k/y_i}{1-y_k/y_i}T_{q,y_i}$$
where $T_{q,y_i}f(q,\st,z_1,\ldots,z_N,y_1,\ldots,y_N):=
f(q,\st,z_1,\ldots,z_N,y_1,\ldots,y_{i-1},qy_i,y_{i+1},\ldots,y_N)$.
It is immediate that after substitution $t=\st/q,\
x_i=y_i/y_{i+1}$ we have $J(q,t,z,t^{-1}x_1^{-1},\ldots,t^{-1}x_{N-1}^{-1})=
f_N(q,\st,z,y)$. It follows that
$\sD'J(q,t,z,x)=(z_1+\ldots+z_N)J(q,t,z,x)$ where
$$\sD':=\sum_{i=1}^Nz_i\prod_{j<i}
\frac{1-q^{-1}t^{i-j-1}x_j\cdots x_{i-1}}{1-t^{i-j}x_j\cdots x_{i-1}}
\prod_{k>i}\frac{1-qt^{k-i+1}x_i\cdots x_{k-1}}{1-t^{k-i}x_i\cdots x_{k-1}}
T_{i,q^{-1}}$$
The theorem follows.
\epr

\ssec{dhl}{Difference equation on $H_\chlam$}
For a weight $\chlam=\sum_{i=1}^{N-1}l_i\check\omega_i$, and $1\leq k\leq N$,
we define $T_k\chlam$ as follows: $T_1\chlam=(l_1-1)\check\omega_1+
l_2\check\omega_2+\ldots+l_{N-1}\check\omega_{N-1},\
T_2\chlam=(l_1+1)\check\omega_1+(l_2-1)\check\omega_2+l_3\check\omega_3+\ldots+
l_{N-1}\check\omega_{N-1},\ \ldots,\ T_N\chlam=l_1\check\omega_1+\ldots+
l_{N-2}\check\omega_{N-2}+(l_{N-1}+1)\check\omega_{N-1}$.
We define the operator $\fD:=\sum_{r=1}^NK_r(\chlam)T_r\chlam$ where
$$K_r(\chlam)=\frac{(1-t^2q^{l_r-1})(1-t^3q^{l_r+l_{r+1}-1})\ldots
(1-t^{N-r+1}q^{l_r+\ldots+l_{N-1}})}
{(1-tq^{l_r})(1-t^2q^{l_r+l_{r+1}})\ldots
(1-t^{N-r}q^{l_r+\ldots+l_{N-1}})}\times$$
$$\times\frac{(1-q^{l_{r-1}+1})(1-tq^{l_{r-1}+l_{r-2}+1})\ldots
(1-t^{r-2}q^{l_{r-1}+\ldots+l_1+1})}{(1-tq^{l_{r-1}})(1-t^2q^{l_{r-1}+l_{r-2}})\ldots
(1-t^{r-1}q^{l_{r-1}+\ldots+l_1})}$$
Now~\refp{brav} and~\reft{shir} admit the following

\cor{diff}
$\fD H_\chlam(q,t,z)=(z_1+\ldots+z_N)H_\chlam(q,t,z)$.
\ecor

\prf
The function
$\fJ(q^{-1},t,wz,q^{\check\lambda})$ on the weight lattice
is an eigenfunction of $\sD$ (with $q$ inverted)
restricted to the weight lattice. According to~\refp{brav},
$H_\chlam(q,t,z)$ is a linear combination of the functions
$\fJ(q^{-1},t,wz,q^{\check\lambda})$ with coefficients independent of $\chlam$.
Hence $H_\chlam(q,t,z)$ is an eigenfunction of $\sD$ (with $q$ inverted)
restricted to the weight lattice (that is $\fD$) as well.
\epr

\ssec{mp}{$H_\chlam$ via Macdonald polynomials}
For $\chlam$ a dominant weight, $P_\chlam$ is the Macdonald
polynomial~\cite{Mac}. To avoid a misunderstanding,
let us state the relation between our
$\chlam=\sum_{i=1}^{N-1}l_i\check\omega_i$, and Macdonald's partitions:
we associate to $(l_1,\ldots,l_{N-1})$ the partition
$(\chlam_N\geq\chlam_{N-1}\geq\ldots\geq\chlam_3\geq\chlam_2\geq0)$
where $l_1=\chlam_2,\ l_2=\chlam_3-\chlam_2,\ \ldots,\ l_{N-1}=\chlam_N-
\chlam_{N-1}$.

\th{h=p}
$$H_\chlam=H_0\prod_{1\leq i\leq j\leq N-1}\frac{(t^{j-i+1};q)_{l_i+\ldots+l_j}}
{(t^{j-i}q;q)_{l_i+\ldots+l_j}}P_\chlam$$
\eth

\prf
The Pieri rule~\cite[Equation~(6.24)(iv) at page~341]{Mac} reads
$\CalD P_\chlam=(z_1+\ldots+z_N)P_\chlam$ where
$\CalD:=\sum_{r=1}^NL_r(\chlam)T_r\chlam$, and
$$L_r(\chlam)=\frac{(1-t^2q^{l_r-1})(1-t^3q^{l_r+l_{r+1}-1})\ldots
(1-t^{N+1-r}q^{l_r+\ldots+l_{N-1}-1})}
{(1-tq^{l_r-1})(1-t^2q^{l_r+l_{r+1}-1})\ldots
(1-t^{N-r}q^{l_r+\ldots+l_{N-1}-1})}\times$$
$$\times\frac{(1-q^{l_r})(1-tq^{l_r+l_{r-1}})\ldots
(1-t^{N-1-r}q^{l_r+\ldots+l_{N-1}})}{(1-tq^{l_r})(1-t^2q^{l_r+l_{r-1}})\ldots
(1-t^{N-r}q^{l_r+\ldots+l_{N-1}})}$$
Since $H_\chlam(q,t,z)$ is an eigenfunction of $\fD$, the function
$$P'_\chlam(q,t,z):=H_\chlam(q,t,z)H_0^{-1}\prod_{1\leq i\leq j\leq N-1}
\frac{(t^{j-i+1};q)_{l_i+\ldots+l_j}^{-1}}{(t^{j-i}q;q)_{l_i+\ldots+l_j}^{-1}}$$
is an eigenfunction of $\CalD$. It vanishes outside the cone of dominant
weights according to~\reft{vanis}, and it equals 1 at $\chlam=0$.
These properties uniquely characterize the Macdonald polynomials
$P_\chlam(q,t,z)$.
\epr

\rem{Weyl}
In view of~\reft{h=p},~\refp{brav} expressing the Macdonald polynomials
in terms of the Baker-Akhiezer function $\fJ(q,t,z,x)$ is nothing but the
generalized Weyl formula~\cite[Proposition~5.3]{ES},~\cite[Theorem~3.9]{CE}.
\erem

\sec{spec}{Speculations for arbitrary simple groups}

\ssec{pcs}{Perverse coherent sheaves}
Let $G$ be an almost simple simply connected complex group with Lie algebra
$\fg$.
We will follow the notations of~\cite{BF11},~\cite{BF12}.
Recall the infinite type scheme $_\fg\bQ$ introduced
in~\cite[Section~2.2]{BF12}:
the quotient by the action of the Cartan torus $T\subset G$ of the space of
maps from $\on{Spec}R=\on{Spec}\BC[[{\mathbf t}^{-1}]]$ to the affinization
of the base affine space $\overline{G/U_-}$ taking value in $G/U_-$ at the
generic point. It is equipped with the action of the proalgebraic group
$G(R)$; the open orbit $_\fg\bQ_\infty=\ {}_\fg\bQ^0$ is nothing but
$G(R)/T\cdot U_-(R)$: the maps taking value in $G/U_-$ at the closed point
$r\in\on{Spec}R$. We denote by $\fj$ the open embedding of $_\fg\bQ^0$
into $_\fg\bQ$.
All the $G(R)$-orbits in $_\fg\bQ$ are numbered by the
defects at $r$ taking value in the cone of positive coroots
$\Lambda_+$ of $G:\ {}_\fg\bQ=\bigsqcup_{\alpha\in\Lambda_+}\ {}_\fg\bQ^\alpha$.
The codimension of $_\fg\bQ^\alpha$ in $_\fg\bQ$ equals $2|\alpha|$.

We introduce the perversity $p(\ {}_\fg\bQ^\alpha)=|\alpha|$; it is immediate
that the function $p$ is strictly monotone and comonotone in the sense
of~\cite{AB}. For a locally free $G(R)\rtimes\BG_m$-equivariant sheaf
$\CF$ on $_\fg\bQ^0$ the construction of~\cite[Section~4]{AB} produces an
object $\fj_{!*}\CF$ of $G(R)\rtimes\BG_m$-equivariant quasicoherent derived
category on $_\fg\bQ$.

%The moduli space
%$\QM^\alpha_\fg$ of degree $\alpha$ quasimaps from $\bC$ to the flag
%variety $\CB_\fg$ is stratified according to the type of defect, see
%e.g.~\cite[Section~1.3]{Ku}. For a stratum
%$\fD_{\beta,\Gamma}\subset\QM^\alpha_\fg$
%we set $p(\eta_{\beta,\Gamma})=-|\beta|-m$ where $m$ is the number of parts of
%the Kostant partition $\Gamma$, and $\eta_{\beta,\Gamma}$ is the generic point
%of $\fD_{\beta,\Gamma}$. If $\eta$ is the generic point of an arbitrary
%reduced irreducible closed subscheme $\fD\subset\QM^\alpha_\fg$, then we set
%$p(\eta)=p(\eta_{\beta,\Gamma})$ where $\fD_{\beta,\Gamma}$ is the
%minimal stratum
%whose closure contains $\fD$. It is immediate that the function $p$ on the
%topological space of the scheme $\QM^\alpha_\fg$ is a monotone and comonotone
%perversity in the sense of~\cite{AB}. It is also clear that for a stratum
%$\fD_{\beta,\Gamma}$ lying in the closed subset
%$^{\on{sing}}\QM^\alpha_\fg\subset\QM^\alpha_\fg$ of singular points
%we have $p(\eta_{\beta,\Gamma})>-|\alpha|<-|\beta|=
%\overline{p}(\eta_{\beta,\Gamma})$ (the dual perversity). Let us denote by
%$\fj:\ ^{\on{reg}}\QM^\alpha_\fg\hookrightarrow\QM^\alpha_\fg$ the
%open embedding of the subset formed by the regular points. Then for
%a locally free sheaf $\CF$ on $^{\on{reg}}\QM^\alpha_\fg$ the construction
%of~\cite[Section~4]{AB} produces an object of coherent derived category
%$\IC(\CF)=\fj_{!*}(\CF[|\alpha|])[-|\alpha|]$. Moreover, $\IC(\CF)$ is
%$G\times\BG_m$-equivariant in case $\CF$ is.

\ssec{LR}{Laumon resolution}
In case $G=\SL(N)$ we denote $_\fg\bQ$ by $\bQ$, and we have
a resolution of singularities $\pi:\ \widetilde\bQ\to\bQ$ where
$\widetilde\bQ$ is the moduli space of flags
$0\subset V_1\subset V_2\ldots\subset V_{N-1}\subset R^N$ of free
$R$-modules, $\on{rk}V_i=i$, along with generators of rank 1 free $R$-modules
$v_i\in\Lambda^iV_i$ defined up to multiplication by a scalar in $\BC$. The smoothness of $\widetilde\bQ$ follows from
the equality $\widetilde\bQ\simeq\left(\prod_{i=1}^{N-1}\on{Hom}_{\on{inj}}
(R^i,R^{i+1})\right)/\prod_{i=1}^{N-1}GL^c(i,R)$ where
$\on{Hom}_{\on{inj}}(R^i,R^{i+1})$ stands
for the open subscheme in the scheme (pro- finite dimensional vector space)
$\on{Hom}_R(R^i,R^{i+1})\simeq R^{i(i+1)}$ formed by all the injective
homomorphisms, while $GL^c(i,R)$ stands for the group of $i\times i$ matrices
with coefficients in $R$, and with {\em constant} nonvanishing determinant.
For a point $\phi\in \bQ^\alpha$
we have $\dim\pi^{-1}(\phi)\leq|\alpha|-1$,
see~\cite[Lemma~2.4.6]{Ku}, i.e. the
morphism $\pi$ is {\em very small}.
Hence $\fj_{!*}(\Omega^j_{\bQ^0})=R\pi_*\Omega^j_{\widetilde\bQ}$.

\ssec{eucha}{Euler characteristics for $\bQ$}
Similarly to~\refss{gen}, we consider
the generating function
$[H^\bullet(\bQ,\fj_{!*}(\Omega^\bullet_{\bQ^0})
\otimes\CO(\chlam))]:=\sum_{i,j}(-1)^{i+j}t^j
[H^i(\bQ,\fj_{!*}(\Omega^j_{\bQ^0})
\otimes\CO(\chlam))]=\sum_{i,j}(-1)^{i+j}t^j[H^i(\widetilde\bQ,
\Omega^j_{\widetilde\bQ}\otimes\pi^*\CO(\chlam))]$.

\prop{chi bQ}
For $\chlam=\sum_{i=1}^{N-1}l_i\check\omega_i$ we have
$$[H^\bullet(\widetilde\bQ,\Omega^\bullet_{\widetilde\bQ}\otimes\pi^*\CO(\chlam))]=
H_0\prod_{1\leq i\leq j\leq N-1}\frac{(t^{j-i+1};q)_{l_i+\ldots+l_j}}
{(t^{j-i}q;q)_{l_i+\ldots+l_j}}\prod_{i=1}^{N-2}\left(\frac{(t^i;q)_\infty}
{(qt^{i+1};q)_\infty}\right)^{N-i-1}P_\chlam.$$
\eprop

\prf
Applying the Atiyah-Bott-Lefschetz localization to the fixed points of
$\BG_m\times T$ in $\widetilde\bQ$ we obtain
$$[H^\bullet(\widetilde\bQ,\Omega^\bullet_{\widetilde\bQ}\otimes\pi^*\CO(\chlam))]=
\sum_{w\in W}
z^{w\check\lambda}\left(\sum_{(\theta_{ij})\in\BN^{\sM^{(N)}}}
C_{(\theta_{ij})}(q^{-1},t,wz)q^{\sum_{(i,j)\in\sM^{(N)}}
(l_i+\ldots+l_{j-1})\theta_{ij}}\right)\times$$
$$\prod_{1\leq i<j\leq N}
{(qt z_{w(j)}/z_{w(i)};q)_{\infty } \over (q z_{w(j)}/z_{w(i)};q)_{\infty}}\times
\left({(qt;q)_{\infty } \over (q;q)_{\infty}}\right)^{N-1}\times
\prod_{\check\alpha\in\check{R}{}^+}\frac{1-twz^{\check\alpha}}
{1-wz^{\check\alpha}}$$
It remains to compare~\refe{60} and the formula~\refe{middle} for $H_\chlam$
with the above formula, taking into account~\reft{h=p} and~\reft{junichi}.
\epr

\ssec{euchar}{Euler characteristics for $_\fg\widehat\bQ$}
In case $G$ is of type $BCFG$,
following~\cite[Section~8.2]{BF11}, we consider a simply connected
simply laced group $G'$
with Lie algebra $\fg'$ and its outer automorphism $\sigma$ such that
$\fg=(\fg')^\sigma$ (i.e. $\fg$ is obtained by folding of $\fg'$).
We define the scheme
$_\fg\widehat\bQ$ as a unique irreducible component of the fixed point
subscheme of the automorphism
$\varsigma$ of $_{\fg'}\bQ$ having nonempty intersection with
$_{\fg'}\bQ^0$ (notations of {\em loc. cit.}).
The orbits of $(G'[[{\mathbf t}^{-1}]])^\varsigma$ on
$_\fg\widehat\bQ$ are numbered by $\Lambda_+(\fg)$, and the minimal extension
from $_\fg\widehat\bQ{}^0$ is defined as in~\refss{pcs}. In order to
unify the notation, in the ADE case let us denote $_\fg\bQ$ by
$_\fg\widehat\bQ$ as well.

%\conj{limi}
%For a $G$-weight $\chlam$, the generating functions
%$[H^\bullet(\widehat{\mathcal{QM}}{}^\alpha_\fg,
%\IC(\Omega^\bullet_{^{\on{reg}}\widehat{\mathcal{QM}}{}^\alpha_\fg})
%\otimes\CO(\chlam))]$ tend to a limit $H_\chlam(q,t,z)$ as
%$\alpha\to\infty$.
%\econj

Similarly to~\refp{h0}, we define $H_0(t)$ as
$W(t)\cdot F(t)$ where $W(t)$ is the Poincar\'e polynomial of the flag
variety $\CB_{\check\fg}$, and
$F(t)=\sum F_it^i$ where $F_i$ is the number of unordered
collections of positive roots $\alpha_1,\ldots,\alpha_k,\beta_1,\ldots,\beta_m
\in R^+(\check\fg)$ such that none of $\alpha_1,\ldots,\alpha_k$ is
simple, and $\sum_{j=1}^k(|\alpha_j|-1)+\sum_{l=1}^m(|\beta_l|+1)=i$.
Here $|\beta|:=\langle\beta,\check\rho\rangle$.

\conj{76}
(a) For a nondominant $G$-weight $\chlam$ we have
$[H^\bullet(_\fg\widehat\bQ,\fj_{!*}(\Omega^\bullet_{_\fg\widehat\bQ{}^0})
\otimes\CO(\chlam))]=0$.

(b) For a dominant $G$-weight $\chlam$ we have
$$[H^\bullet(_\fg\widehat\bQ,\fj_{!*}(\Omega^\bullet_{_\fg\widehat\bQ{}^0})
\otimes\CO(\chlam))]=H_0\prod_{\alpha\in R^+(\check\fg)}
\frac{(t^{|\alpha|};q)_{\langle\alpha,\chlam\rangle}}
{(t^{|\alpha|-1}q;q)_{\langle\alpha,\chlam\rangle}}
\prod\frac{(t^{|\alpha|-1};q)_\infty}{(qt^{|\alpha|};q)_\infty}
P_\chlam$$
where $P_\chlam(q,t,z)$ is the Macdonald polynomial for $G$,
and the second product is taken over all {\em nonsimple} positive roots
of $R^+(\check\fg)$.
\econj

\bigskip
\footnotesize{
{\bf A.B.}: Department of Mathematics, Brown University,
151 Thayer St., Providence RI
02912, USA;\\
{\tt braval@math.brown.edu}}

\footnotesize{
{\bf M.F.}: IMU, IITP and National Research University
Higher School of Economics\\
Department of Mathematics, 20 Myasnitskaya st, Moscow 101000 Russia;\\
{\tt fnklberg@gmail.com}}

\footnotesize{
{\bf J.S.}: Graduate School of Mathematical Sciences, University of Tokyo,
Komaba, Tokyo 153-8914, Japan;\\
{\tt shiraish@ms.u-tokyo.ac.jp}}

\end{document}